\documentclass[11pt, notitlepage]{article}
\usepackage{amssymb,amsmath,comment}
\catcode`\@=11 \@addtoreset{equation}{section}
\def\thesection{\arabic{section}}

\def\theequation{\thesection.\arabic{equation}}
\catcode`\@=12
\usepackage{colortbl}
\usepackage{a4wide}
\usepackage{parskip}
\usepackage{xcolor}

\newcommand{\e}{\epsilon}
\newcommand{\pa} {\partial}

\newcommand{\de} {\delta}

\newcommand{\Om} {\Omega}

\newcommand{\De} {\Delta}
\newcommand{\la} {\lambda}
\newcommand{\La} {\Lambda}
\newcommand{\noi} {\noindent}

\newcommand{\mb} {\mathbb}
\newcommand{\mc} {\mathcal}

\usepackage[all]{xy}
\catcode`\@=11
\def\theequation{\@arabic{\c@section}.\@arabic{\c@equation}}
\catcode`\@=12

\newtheorem{Theorem}{Theorem}[section]
\newtheorem{Lemma}[Theorem]{Lemma}
\newtheorem{Proposition}[Theorem]{Proposition}

\newtheorem{Remark}[Theorem]{Remark}
\newtheorem{Definition}[Theorem]{Definition}

\begin{document}
\title{Critical growth elliptic problems with Choquard type nonlinearity}
\author{ Tuhina Mukherjee\footnote{tulimukh@gmail.com}\\Tata Institute of Fundamental Research(TIFR)\\ Centre of Applicable Mathematics,\\ Bangalore, India.\\ K. Sreenadh\footnote{sreenadh@maths.iitd.ac.in}\\Department of Mathematics\\ Indian Institute of Technology Delhi, \\Hauz Khas, New Delhi-110016, India}
\date{}
\maketitle

\abstract{This article deals with a survey of recent developments and results  on Choquard equations where we focus on the existence and multiplicity of solutions of the partial differential equations which involve the nonlinearity of convolution type. Because of its nature, these equations are categorized under the nonlocal problems. We give a brief survey on the work already done in this regard following which we illustrate the problems we have addressed. Seeking the help of variational methods and asymptotic estimates, we prove our main results.}

\section{A Brief Survey}
\label{sec:1}
We devote our first section on briefly glimpsing the results that have already been proved in the context of existence and multiplicity of solutions of the Choquard equations. Consider the problem
\begin{equation}\label{BS-1}
-\Delta u +u = (I_\alpha \ast |u|^p)|u|^{p-2}u \; \text{in} \; \mathbb{R}^n
\end{equation}
where $u: \mathbb R^n \to \mathbb{R}$ and $I_\alpha: \mathbb R^n \to \mathbb R$ is the {Riesz} potential defined by
\[I_\alpha(x) = \frac{\Gamma\left(\frac{n-\alpha}{2}\right)}{\Gamma\left(\frac{\alpha}{2}\right)\pi^{\frac{n}{2}} 2^\alpha |x|^{n-\alpha}}\]
for $\alpha \in (0,n)$ and $\Gamma$ denotes the Gamma function. Equation \eqref{BS-1} is generally termed as Choquard equations or the Hartree type equations. It has various physical significance. In the case $n=3$, $p=2$ and $\alpha=2$,  \eqref{BS-1} finds its origin in a work by S.I. Pekar describing the quantum
mechanics of a polaron at rest  \cite{Pekar}. Under the same assumptions, in $1976$ P. Choquard used \eqref{BS-1} to describe an electron trapped in its own hole, in a certain approximation to Hartree-Fock theory of one component plasma \cite{Choquard}. Following standard critical point theory, we expect that solutions of \eqref{BS-1} can be viewed as critical points of the energy functional
\[J(u)= \frac12\int_{\Om}(|\nabla u|^2+ u^2)- \frac{1}{2p}\int_\Om(I_\alpha*|u|^p)|u|^p. \]
It is clear from the first term that naturally we have to take $u \in H^1(\mb R^n)$ which makes the first and second term well defined. Now the question is whether the third term is well defined and sufficiently smooth over $H^1(\mathbb{R}^n)?$ For this, we recall the following Hardy-Littlewood- Sobolev inequality.
\begin{Theorem}\label{HLS} Let $t,r>1$ and $0<\mu<n $ with $1/t+\mu/n+1/r=2$, $f \in L^t(\mathbb R^n)$ and $h \in L^r(\mathbb R^n)$. Then there exists a constant $C(t,n,\mu,r)$, independent of $f,h$ such that
	\begin{equation}\label{har-lit}
	\int_{\mb R^n}\int_{\mb R^n} \frac{f(x)h(y)}{|x-y|^{\mu}}\mathrm{d}x\mathrm{d}y \leq C(t,n,\mu,r)\|f\|_{L^t}\|h\|_{L^r}.
	\end{equation}
	{ If $t =r = \textstyle\frac{2n}{2n-\mu}$ then
		\[C(t,n,\mu,r)= C(n,\mu)= \pi^{\frac{\mu}{2}} \frac{\Gamma\left(\frac{n}{2}-\frac{\mu}{2}\right)}{\Gamma\left(n-\frac{\mu}{2}\right)} \left\{ \frac{\Gamma\left(\frac{n}{2}\right)}{\Gamma(n)} \right\}^{-1+\frac{\mu}{n}}.  \]
		The inequality in \eqref{har-lit} is achieved  if and only if $f\equiv (constant)h$ and
		\[h(x)= A(\gamma^2+ |x-a|^2)^{\frac{-(2n-\mu)}{2}}\]
		for some $A \in \mathbb C$, $0 \neq \gamma \in \mathbb R$ and $a \in \mathbb R^n$.}
\end{Theorem}
For $u \in H^1(\mb R^n),$ let $f = h= |u|^p$,  then by Theorem \ref{HLS},
\[ \int_{\mb R^n}\int_{\mb R^n} \frac{|u(x)|^p|u(y)|^p}{|x-y|^{n-\alpha}}\mathrm{d}x\mathrm{d}y \leq C(t,n,\mu,p)\left(\int_{\mb R^n}|u|^{\frac{2np}{n+\alpha}}\right)^{1+\frac{\alpha}{n}}. \]
This is well defined if $u \in L^{\frac{2np}{n+\alpha}}(\mb R^n)$. By the classical Sobolev embedding theorem, the embedding $H^1(\mb R^n) \hookrightarrow L^r(\mb R^n)$ is continuous when $r \in [2, 2^*)$, where $2^*= \frac{2n}{n-2}$. This implies $u \in L^{\frac{2np}{n+\alpha}}(\mb R^n)$ if and only if
\begin{equation}\label{cond-on-p}
2_\alpha:=\frac{n+\alpha}{n}\leq p \leq \frac{n+\alpha}{n-2}:= 2^*_\alpha.
\end{equation}
The constant $2_\alpha$ is termed as the lower critical exponent and $2_\alpha^*$ is termed as the upper critical exponent in the sense of Hardy-Littlewood-Sobolev inequality.
Then we have the following result.
\begin{Theorem}
	If $p \in (1,\infty)$ satisfies \eqref{cond-on-p}, then the functional $J$ is well-defined and continuously Fr\'{e}chet differentiable
	on the Sobolev space $H^1(\mb R^n)$.
	Moreover, if $p\geq 2$, then the functional $J$ is twice continuously Fr\'{e}chet differentiable.
\end{Theorem}
This suggests that it makes sense to define the solutions of \eqref{BS-1} as critical points of $J$. A remarkable feature in Choquard nonlinearity is the appearance of a lower nonlinear restriction: the lower critical exponent $2_\alpha >1$. That is the nonlinearity is superlinear.
\subsection{Existence and multiplicity results}
\begin{Definition}
	A function $u \in H^1(\mb R^n)$ is said to be a weak solution of \eqref{BS-1} if it satisfies
	\[\int_{\mb R^n }(\nabla u \nabla v + uv)~dx + \int_{\mb R^n} \left(\int_{\mb R^n}\frac{|u(y)|^p}{|x-y|^{\alpha}}~dy \right)|u|^{p-2}uv ~dx=0\]
	for each $v \in H^1(\mb R^n)$.
\end{Definition}
\begin{Definition}
	We define a solution $u \in H^1(\mb R^n)$ to be a
	groundstate of the Choquard equation  \eqref{BS-1} whenever it is a solution that
	minimizes the functional $J$ among all nontrivial solutions.
\end{Definition}
In \cite{vs} V. Moroz and J. Van Schaftingen studied the existence of groundstate solutions and their asymptotic behavior  using concentration-compactness lemma. The groundstate solution has been identified as infimum of $J$ on the Nehari manifold
\[\mc N =\{u \in H^1(\mb R^n): \; \langle J^\prime(u),u\rangle =0\}\]
which is equivalent to prove that the mountain pass minimax level
$\displaystyle \inf_{\gamma \in \Gamma}\sup_{[0,1]}J\circ\gamma$ is a critical value. Here
the class of paths $\Gamma$ is defined by
$\Gamma = \{\gamma \in C([0,1];H^1(\mb R^n)):\; \gamma(0)=0 \text{ and } J(\gamma(1))<0\}$.
Precisely, they proved the following existence result-
\begin{Theorem}
	If $2_\alpha<p<2^*_\alpha$ then there exists a nonzero weak solution $u \in W^{1,2}(\mb R^n)$ of \eqref{BS-1} which is a groundstate solution of \eqref{BS-1}.
\end{Theorem}
They have also proved the following Pohozaev identity:
\begin{Proposition}\label{vs-poho}
	Let $u \in H^{2}_{loc}(\mb R^n)\cap W^{1,\frac{2np}{n+\alpha}}{(\mb R^n)}$ is a weak solution of the equation
	\[-\De u + u = (I_\alpha*|u|^p)|u|^{p-2}u\; \text{in}\; \mathbb{R}^n\]
	then
	\[\frac{n-2}{2}\int_{\mb R^n}|\nabla u|^2+\frac{n}{2}\int_{\mb R^n}|u|^2= \frac{n+\alpha}{2p}\int_{\mb R^n}(I_\alpha*|u|^p)|u|^p.\]
\end{Proposition}
Pohozaev identity for some Choquard type nonlinear equations has also been studied in \cite{menzala,css}. Using Proposition \ref{vs-poho}, they proved the following nonexistence result-
\begin{Theorem}\label{nonext}
	If $p \leq 2_\alpha$ or $p\geq 2^*_\alpha$ and  $u \in H^{1}(\mb R^n)\cap L^{\frac{2np}{n+\alpha}}(\mb R^n)$ such that $\nabla u \in  H^{1}_{loc}(\mb R^n)\cap L^{\frac{2np}{n+\alpha}}_{loc}(\mb R^n)$  satisfies \eqref{BS-1} weakly, then  $u \equiv 0$.
\end{Theorem}
Next important thing to note is the following counterpart of Brezis-Leib lemma: \\ If the sequence $\{u_k\}$ converges weakly to $u$ in $H^1(\mb R^n)$, then
\begin{equation}\label{BS-2}
\lim_{k \to \infty} \int_{\mb R^n}(I_\alpha*|u_k|^p)|u_k|^p- (I_\alpha*|u-u_k|^p)|u-u_k|^p=  \int_{\mb R^n}(I_\alpha*|u|^p)|u|^p.
\end{equation}
One can find its proof in \cite{mms,vs}. Equation \eqref{BS-2} plays a crucial role in obtaining the solution where there is a lack of compactness.\\
Next coming to the positive solutions, in \cite{AFY} authors studied the existence of solutions for the following  equation
\begin{equation}\label{BS-3}
-\Delta u +V(x)u = (|x|^{-\mu}\ast F(u))f(u), \; u>0\; \text{in}\; \mb R^n, \; u \in D^{1,2}(\mb R^n)
\end{equation}
where $F$ denotes primitive of $f$, $n\geq 3$ and $\mu \in (0,n)$. Assumptions on the potential function $V$ and the function $f$ are as follows:
\begin{enumerate}
	\item[(i)] $\displaystyle \lim_{s \to 0^+} \frac{sf(s)}{s^q}<+\infty$ for $q \geq 2^* =\frac{2n}{n-2}$
	\item[(ii)] $\displaystyle \lim_{s \to \infty} \frac{sf(s)}{s^p}=0$ for some $p \in \left( 1, \frac{2(n-\mu)}{n-2}\right)$ when $\mu \in (1,\frac{n+2}{2})$,
	\item[(iii)] There exists $\theta>2$ such that $1<\theta F(s)<2f(s)s$ for all $s>0$,
	\item[(iv)] $V$ is a nonnegative continuous function.
\end{enumerate}
Define the function $\mc V : [1,+\infty) \to [0,\infty)$ as
\[ \mc V(R) = \frac{1}{R^{(q-2)(n-2)}}\inf_{|x|\geq R} |x|^{(q-2)(n-2)}V(x)\]
Motivated by the articles \cite{BL,BGM1,BGM2}, authors proved the following result in \cite{AFY}.
\begin{Theorem}
	Assume that $0<\mu< \frac{n+2}{2}$ and $(i)-(iv)$ hold. If there exists a constant $\mc V_0>0$ such that if $\mc V(R)>\mc V_0$ for some $R>1$, then \eqref{BS-3} admits a positive solution.
\end{Theorem}
Taking $p=2$ in \eqref{BS-1}, Ghimenti, Moroz and Schaftingen \cite{nodal}
established existence of a least action nodal solution which appeared as the limit of  minimal action nodal solutions for \eqref{BS-1} when $p \searrow 2$. They proved the following theorem  by constructing a Nehari nodal set and minimizing the corresponding energy functional over this set.
\begin{Theorem}
	If $\alpha \in ((n-4)^+, n)$ and $p=2$ then \eqref{BS-1} admits a least action nodal solution.
\end{Theorem}
In \cite{ZHZ}, Zhang et al. proved the existence of infinitely many distinct solutions for the following generalized Choquard equation using the index theory
\begin{equation}\label{BS-5}
-\Delta u +V(x)u = \left( \int_{\mb R^n} \frac{Q(y)F(u(y))}{|x-y|^{\mu}}~dy\right)Q(x)f(u(x)) \; \text{in}\; \mb R^n
\end{equation}
where $\mu \in (0, n)$, $V$ is periodic, $f$ is either odd or even and some additional assumptions.
Although Theorem \ref{nonext} holds, when $p=2_\alpha$ in \eqref{BS-1}, Moroz and Schaftingen in \cite{vs2} proved some existence and nonexistence of solutions for the problem
\begin{equation}\label{BS-4}
-\Delta u + V(x)u = (I_\alpha \ast |u|^{2_\alpha})|u|^{2_\alpha-2}u \; \text{in} \; \mathbb{R}^n
\end{equation}
where the potential $V \in L^\infty(\mb R^n)$ and must not be a constant. They proved existence of a nontrivial solution if \[\liminf_{|x|\to \infty}(1-V(x))|x|^2 >  \frac{n^2(n-2)}{4(n+1)}\] and gave necessary conditions for existence of solutions of \eqref{BS-4}. Because $2_\alpha$ is the lower critical exponent in the sense of Theorem \ref{HLS}, a lack of compactness occurs in minimization technique. So concentration compactness lemma and Brezis Lieb type lemma plays an important role. Equation \eqref{BS-4} was reconsidered by Cassani, Schaftingen and Zhang in \cite{CSZ} where they gave necessary and sufficient condition for the existence of positive ground state solution depending on the potential $V$.\\
Very recently, in \cite{divsree}, authors studied some existence and multiplicity results for the following {critical} growth Kirchhoff- Choquard equations
\[ -M(\|u\|^2)\Delta u = \la u + (I_\alpha \ast |u|^{2_{\mu}^{*}})|u|^{2_{\mu}^{*}-2}u \; \text{in} \; \Om\textcolor{red}{,}\; u=0 \; \text{on}\; \pa \Om\]
where $M(t)\sim at+bt^\theta, \theta\ge 1$ for some constants $a$ and $b$.

\noi Now let us consider the critical dimension case that is $n=2$ commonly known as the Trudinger-Moser case. When $n=2$, the critical Sobolev exponent becomes infinity and the embedding goes as $W^{1,2}(\mb R^2) \hookrightarrow L^q(\mb R^2)$ for $q \in [2,\infty)$ whereas $W^{1,2}(\mb R^2)\not\hookrightarrow L^\infty(\mb R^2)$. The following \textit{Trudinger-Moser inequality} plays a crucial role when $n=2$.
\begin{Theorem}\label{TM-ineq}
	For $u \in W^{1,2}_0(\mb R^2)$,
	\[\int_{\mb R^2} [\exp(\alpha|u|^{2})-1]~dx < \infty.\]
	Moreover if $\|\nabla u\|_2\leq 1$, $\|u\|_2 \leq M$ and $\alpha < 4\pi$ then there exists a $C(\alpha,M)>0$ such that
	\[\int_{\mb R^2} [\exp(\alpha|u|^{2})-1]~dx <C(M,\alpha).\]
\end{Theorem}
Motivated by this, the nonlinearity in this case is an appropriate exponential function. The following singularly perturbed Choquard equation
\begin{equation}\label{BS-8}
-\epsilon^2 \Delta u +V(x)u = \e^{\mu-2} \left( |x|^{-\mu}\ast F(u)\right)f(u) \;\text{in}\; \mb R^2
\end{equation}
was studied by Alves et al. in \cite{yang-JDE}. Here $\mu\in (0,2)$, $V$ is a continuous potential, $\e$ is a positive parameter, $f$ has critical exponential growth in the sense of Trudinger-Moser and $F$ denotes its primitive. Under appropriate growth assumptions on $f$, authors in \cite{yang-JDE} proved existence of a ground state solution to \eqref{BS-8} when $\e=1$ and $V$ is periodic and also established the existence and concentration of semiclassical ground state solutions of \eqref{BS-8} with respect to $\e$. An existence result for Choquard equation with exponential nonlinearity in $\mb R^2$ has also been proved in \cite{yang-JCA}. The Kirchhoff-Choquard problems in this case are studied in the work \cite{arora}.\\

\noi Now let us consider the Choquard equations in the bounded domains. In particular, consider the Brezis-Nirenberg type problem for Choquard equation
\begin{equation}\label{BS-6}
-\De u = \la u + \left( \int_{\Om}\frac{|u|^{2^*_\mu}(y)}{|x-y|^\mu}~dy\right) |u|^{2^*_\mu-2}u\; \text{in}\; \Om, \; u = 0 \; \text{on}\; \partial\Om
\end{equation}
where $\Om$ is bounded domain in $\mb R^n$ with Lipschitz boundary, $\la \in \mb R$ and $2^*_\mu = \displaystyle\frac{2n-\mu}{n-2}$ which is the critical exponent in the sense of Hardy-Littlewood-Sobolev inequality. These kind of problems are motivated by the celebrated paper of Brezis and Nirenberg \cite{breniren}. Gao and Yang in \cite{gy1} proved existence of nontrivial solution to \eqref{BS-6} for $n \geq 4$ in case $\la$ is  not an eigenvalue of $-\De$ with Dirichlet boundary condition and for a suitable range of $\la$ when $n=3$. They also proved the nonexistence result when $\Om$ is a star shaped region with respect to origin. Here,  the best constant for the embedding is defined as
\[S_{H,L}:= \inf\left\{\int_{\mb R^n}|\nabla u|^2: \; u\in H^1(\mb R^n), \; \int_{\mb R^n}(|x|^{-\mu}*|u|^{2^*_\mu})|u|^{2^*_\mu}dx=1\right\}.\]
They showed that the  minimizers of $S_{H,L}$ are of the form $U(x)= \left( \frac{b}{b^2+|x-a|^2}\right)^{\frac{n-2}{2}}$ where $a,b$ are appropriate constants. We remark that $U(x)$ is  the Talenti function which also forms minimizers of $S$, the best constant in the embedding $H^1_0(\Om)$ into $L^{2^*}(\Om)$. Let us consider the family $U_{\e}(x)= \epsilon^{\frac{2-n}{2}}U(\frac{x}{\epsilon})$.  Using Brezis-Lieb lemma, in \cite{gy1} it was  shown that every Palais Smale sequence is bounded and the first  critical level is
\[c < \frac{n+2-\mu}{4n-2\mu}S_{H,L}^{\frac{2n-\mu}{n+2-\mu}}.\]
If $Q_\la:= \inf\limits_{u \in H^1_0(\Om)\setminus \{0\}}\frac{\int_{\Om}|\nabla u|^2- \la u}{\int_{\Om R^n}(|x|^{-\mu}*|u|^{2^*_\mu})|u|^{2^*_\mu}dx}$, then $Q_\la < S_{H,L}$ which can be shown using $U_\epsilon$'s. Then using Mountain pass lemma and Linking theorem depending on the dimension $n$, existence of first solution to \eqref{BS-6} is shown. The nonexistence result was proved after establishing a Pohozaev type identity.
Gao and Yang also studied Choquard equations with concave-convex power nonlinearities in \cite{gy2} with Dirichlet boundary condition.\\
Very recently, the effect of topology of domain on the solution of Choquard equations has been studied by some researchers. Ghimenti and Pagliardini \cite{GP} proved that the number of positive solution of the following Choquard equation
\begin{equation}\label{BS-7}
-\Delta u -\la u= \left( \int_\Om \frac{|u|^{p_\e}(y)}{|x-y|^{\mu}}~dy\right)|u|^{p_\e-2}u,\;u>0 \; \text{in}\;\Om, \; \; u=0 \;\text{in}\; \mb R^n \setminus \Om
\end{equation}
depends on the topology of the domain when the exponent $p_\e$ is very close to the critical one. Precisely, they proved-
\begin{Theorem}
	There exists $\bar \e > 0$ such that for every $\e \in (0, \bar \e]$, Problem \eqref{BS-7} has at least $cat_\Om(\Om)$ low energy solutions. Moreover, if
	it is not contractible, then there exists another solution with higher energy.
\end{Theorem}
Here $cat_\Om(\Om)$ denotes the Lusternik-Schnirelmann category of $\Om$. They used variational methods to look for critical points of a suitable functional and proved a multiplicity result through category methods. This type of result was historically introduced by  Coron
for local problems in \cite{Bahri}. Another significant result in this regard has been recently obtained by authors in \cite{divya}. Here they showed existence of a high energy solution for
\[-\Delta u =  \left( \int_\Om \frac{|u|^{2^*_\mu}(y)}{|x-y|^{\mu}}~dy\right)|u|^{2^*_\mu-2}u\; \text{in}\; \Om,\;u=0\;\text{on}\; \partial \Om,\]
where $\Om$ is an annular type domain with sufficiently small inner hole.

\subsection{Radial symmetry and Regularity of solutions}
Here, we try to give some literature on radially symmetric solutions and regularity of weak solutions constructed variationally for Choquard equations.

\noi First we come to the question of radially symmetric solutions. Is all the positive solutions for the equation
\begin{equation}
\Delta u -\omega u + (|x|^{-\mu} \ast |u|^{2_\alpha})p|u|^{2_\alpha-2}u=0, \; \omega >0,\; u \in H^1(\mb R^n)
\end{equation}
are radially symmetric and monotone decreasing about some fixed point? This was an open problem which was settled by Ma and Zhao \cite{MZ} in case $2\leq p < \frac{2n-\mu}{n-2}$ and some additional assumptions.  The radial symmetry and uniqueness of minimizers corresponding to some Hartree equation has also been investigated in \cite{JV}. Recently, Wang and Yi \cite{WY}  proved that if $u \in C^2(\mb R^n)\cap H^1(\mb R^n)$ is a positive radial solution of \eqref{BS-1} with $p=2$ and $\alpha =2$ then $u$ must be unique. Using Ma and Zhao's result, they also concluded that the positive solutions of \eqref{BS-1} in this case is uniquely determined, up to translations in the dimension $n=3,4,5$. Huang et al. in \cite{HYY} proved that \eqref{BS-1} with $n=3$ has at least one radially symmetric solutions changing sign exactly $k$-times for each $k$ when $p \in \left(2.5,5 \right)$. Taking $V \equiv 1$ in \eqref{BS-5} and  $f$ satisfies almost necessary the upper critical growth conditions in the spirit of Berestycki and Lions, Li and Tang \cite{LT} very recently proved that \eqref{BS-5} has a ground state solution, which has constant sign and is radially symmetric with respect to some point in $\mb R^n$. They used the Pohozaev manifold and a compactness lemma by Strauss to conclude their main result. For further results regarding Choquard equations, we suggest readers to refer \cite{survey} which  extensively covers  the existing literature on the topic. Very recently, in \cite{divya} authors studied the classification problem and proved that all positive solutions of the following equation are radially symmetric: for $p=2_{\mu}^{*}$
\begin{equation} \label{co2} -\De u= (I_{\mu}*|u|^{p}) |u|^{p-2} u \; \text{in} \; \mathbb{R}^n. \end{equation}
They observed that the solutions of this problem satisfy the integral system of equations
\begin{equation}
\begin{aligned}
& u(x)=\int_{\mathbb{R}^n}\frac{u^{p-1}(y)v(y)}{|x-y|^{N-2}}~dy,  u\geq 0 \text{ in } \mathbb{R}^n \\
& v(x)=\int_{\mathbb{R}^n}\frac{u^p(y)}{|x-y|^{N-\mu}}~dy,  v\geq 0 \text{ in } \mathbb{R}^n.
\end{aligned}
\end{equation}	
By obtaining the regularity estimates and using moving method they proved the following result:
\begin{Theorem}\label{cothm5}
	Every non-negative solution $u \in  D^{1,2}(\mathbb{R}^N)$ of equation \eqref{co2} is radially symmetric, monotone decreasing and of the form
	\begin{align*}
	u(x)= \left(\frac{c_1}{c_2+|x-x_0|^2}\right)^{\frac{N-2}{2}}
	\end{align*}
	for some constants $c_1,c_2>0$ and some $x_0 \in \mathbb{R}^N$.
\end{Theorem}

\noi Next we recall some regularity results for the problem \eqref{BS-1}. Fix $\alpha \in (0,n)$ and consider the problem \eqref{BS-1}, then in \cite{vs} authors showed the following-
\begin{Theorem}\label{reg-1}
	If $u\in H^1(\mb R^n)$ solves \eqref{BS-1} weakly for $p \in (2_\alpha, 2_\alpha^*)$ then $u \in L^1(\mb R^n) \cap C^2(\mb R^n)$, $u \in W^{2,s}(\mb R^n)$ for every $s>1$ and $u \in C^\infty(\mb R^n \setminus u^{-1}\{0\})$.
\end{Theorem}
The classical bootstrap method for subcritical semilinear elliptic problems
combined with estimates for Riesz potentials allows them to prove this result. Precisely, they first proved that $I_\alpha \ast |u|^p \in L^\infty(\mb R^n)$ and using the Calder\'{o}n-Zygmund theory  they obtain  $u \in W^{2,r}(\mb R^n)$ for every $r>1$. Then the proof of Theorem \ref{reg-1} followed from application of Morrey--Sobolev embedding and  classical Schauder regularity estimates. In \cite{MS-ams}, author extended a special case of the
regularity result by Brezis and Kato \cite{BK} for the Choquard equations. They proved the following-
\begin{Theorem}
	If $H, K \in L^{\frac{2n}{\alpha}}(\mb R^n) \cap L^{\frac{2n}{\alpha+2}}(\mb R^n)$ and $u \in H^1(\mb R^n)$ solves
	\[-\Delta u + u= (I_\alpha \ast Hu)K\; \text{in}\; \mathbb{R}^n\]
	then $u \in L^p(\mb R^n)$ for every $p \in \left[2, \frac{2n^2}{\alpha(n-2)}\right)$.
\end{Theorem}
They proved it by establishing a nonlocal counterpart of Lemma $2.1$ of \cite{BK} in terms of the Riesz potentials. After this, they showed that the convolution term is a bounded function that is $I_\alpha \ast |u|^p \in L^\infty(\mb R^n)$. Therefore,
\[|-\Delta u + u| \leq C (|u|^{\frac{\alpha}{n}}+ |u|^{\frac{\alpha+2}{n-2}})\]
that is the right hand side now has subcritical growth with respect to the Sobolev embedding. So by the classical bootstrap method for subcritical local problems in bounded domains,
it is  deduced that $u \in W^{2,p}_{loc}(\mb R^n)$ for every $p \geq 1$. Moreover it holds that if \eqref{BS-1} admits a positive solution and $p$ is an even integer then $u \in C^\infty$, refer \cite{Lei1,Lei2,MSS1}. Using appropriate test function and  results from \cite{MS-ams}, Gao and Yang in \cite{gy2} established the following regularity and $L^\infty$ estimate for problems on bounded domains-
\begin{Lemma}
	Let $u$ be the solution of the problem
	\begin{equation}
	-\Delta u = g(x,u)\;\text{in}\; \Om,\;\;
	u \in H^{1}_0(\Om),
	\end{equation}
	where $g$ is satisfies $|g(x,u)| \leq C(1+|u|^p)+ \left(\displaystyle\int_\Om \frac{|u|^{2^*_\mu}}{|x-y|^{\mu}}dy\right)|u|^{2^*_\mu-2}u$, $\mu \in (0,n)$, $1<p<2^*-1$ and $C>0$ then $u \in L^\infty(\Om)$.
\end{Lemma}
As a consequence of this lemma we can obtain  $u \in C^2(\bar \Om)$ by adopting the classical $L^p$ regularity theory of elliptic equations.

\subsection{Choquard equations involving the $p(.)$-Laplacian}\label{sec-1.2}
Firstly, let us consider the quasilinear generalization of the Laplace operator that is the $p$-Laplace operator defined as
\[-\De_p u:=-\nabla\cdot(|\nabla u|^{p-2}\nabla u), \; 1<p<\infty.\]
The Choquard equation involving $-\De_p$ has been studied in \cite{ay1,ay2,ay3}. In \cite{ay1}, Alves and Yang studied concentration behavior of solutions for the following quasilinear Choquard equation
\begin{equation}\label{BS-9}
-\epsilon^p \Delta_p u +V(x)|u|^{p-2}u = \e^{\mu-n} \left( \int_{\mb R^n}\frac{Q(y)F(u(y))}{|x-y|^{\mu}}\right)Q(x)f(u) \;\text{in}\; \mb R^n
\end{equation}
where $1<p<n$, $n \geq 3$, $0 < \mu < n$, $V$ and $Q$ are two continuous real valued functions on $\mb R^n$,
$F(s)$ is the primitive function of $f(s)$ and $\e$ is a positive parameter. Taking $Q \equiv 1$, Alves and Yang also studied \eqref{BS-9} in \cite{ay2}.
Recently, Alves and Tavares proved a version of Hardy-Littlewood-Sobolev inequality with variable exponent in \cite{p(x)-choq} in the spirit of variable exponent Lebesgue and Sobolev spaces. Precisely, for $p(x), q(x)\in C^+(\mathbb{R}^n)$ with $p^-:=\min\{p(x),0\}>1,$ and $ q^{-}>1$, the following holds:
\begin{Theorem}\label{var-exp}
	Let  $h \in L^{p^+}(\mb R^n) \cap L^{p^-}(\mb R^n)$, $g \in L^{q^+}(\mb R^n) \cap L^{q^-}(\mb R^n)$ and $\la :{\mb R^n} \times \mb R^n \to \mb R$ be a continuous function such that $0\leq \la^+\leq \la^-<n$ and
	\[\frac{1}{p(x)}+ \frac{\la(x,y)}{n}+ \frac{1}{q(y)}=2,\; \forall x,y\in \mb R^n.\]
	Then there exists a constant $C$ independent of $h$ and $g$ such that
	\[\left|\int_{\mb R^n}\int_{\mb R^n} \frac{h(x)g(y)}{|x-y|^\mu}~dxdy\right| \leq C(\|h\|_{p^+}\|g\|_{q^+}+ \|h\|_{p^-}\|g\|_{q^-}). \]
\end{Theorem}
In the spirit of Theorem \ref{var-exp}, authors in \cite{p(x)-choq} proved existence of a solution $u \in W^{1,p(x)}(\mb R^n)$ to the following quasilinear Choquard equation using variational methods under the subcritical growth conditions on $f(x,u)$:
\begin{equation}\label{BS-12}
-\Delta_{p(x)}u + V(x)|u|^{p(x)-2}u=\left(\frac{F(x,u(x))}{|x-y|^{\la(x,y)}}~dx \right)f(y,u(y))\; \text{in}\; \mb R^n,\;
\end{equation}
where  $\Delta_{p(x)}$ denotes the $p(x)$-Laplacian defined as
$-\De_{p(x)} u:=-div(|\nabla u|^{p(x)-2}\nabla u)$, $V$, $p$, $f$ are real valued continuous functions and $F$ denotes primitive of $f$ with respect to the second variable.
\section{Choquard equations involving the fractional Laplacian}
In this section, we summarize our contributions related to the existence and multiplicity results concerning different Choquard equations, in separate subsections. We employ the variational methods and used some asymptotic estimates to achieve our goal. While dealing with critical exponent in the sense of Hardy-Littlewood-Sobolev inequality, we always consider the upper critical exponent. We denote $\|\cdot\|_r$ as the $L^r(\Om)$ norm.
\subsection{Brezis-Nirenberg type existence results}
The fractional Laplacian operator $(-\De)^s$ on the set of the Schwartz class functions is defined  as
\[ (-\De)^s u(x) = -\mathrm{P.V.}\int_{\mb R^n} \frac{u(x)-u(y)}{\vert x-y\vert^{n+2s}}~{d}y\]
({up to a normalizing constant}), where $\mathrm{P.V.}$ denotes the Cauchy principal value, $s \in (0,1)$ and $n>2s$. The operator $(-\De )^s$  is the infinitesimal generator of L$\acute{e}$vy stable diffusion process. The equations involving this operator arise in the modelling of anomalous diffusion in plasma, population dynamics, geophysical fluid dynamics, flames propagation, chemical reactions in liquids and American options in finance. Motivated by \eqref{BS-6}, in \cite{TS-1} we considered the following doubly nonlocal equation involving the fractional Laplacian with noncompact nonlinearity
\begin{equation}\label{BS-10}
(-\De)^su = \left( \int_{\Om}\frac{|u|^{2^*_{\mu,s}}}{|x-y|^{\mu}}\mathrm{d}y \right)|u|^{2^*_{\mu,s}-2}u +\la u \; \text{ in } \Om, \; \; u =0 \; \text{ in } \mathbb R^n\setminus \Om,
\end{equation}
where $\Om$ is a bounded domain in $\mathbb R^n$ with Lipschitz boundary, $\la $ is a real parameter, $s\in (0,1)$, $2^*_{\mu,s}= \displaystyle\frac{2n-\mu}{n-2s}$, {$0<\mu<n$} and $n>2s$. Here, $2^*_{\mu,s}$ appears as the upper critical exponent in the sense of Hardy-Littlewood-Sobolev inequality when the function is taken in the fractional Sobolev space $H^s(\mb R^n):= \{u \in L^2(\mb R^n): \|(-\De)^{\frac{s}{2}}u\|_2 < \infty\}$ which is continuously embedded in $L^{2^*_s}(\mb R^n)$ where $2^*_s= \displaystyle\frac{2n}{n-2s}$. For more details regarding the fractional Sobolev spaces and its embeddings, we refer \cite{hitch}. Following are the main results that we have proved-
\begin{Theorem}\label{thrm1}
	{Let  $n \geq 4s$ for $s \in (0,1)$, then \eqref{BS-10} has a positive weak solution for every $\la>0$ such that $\la$ is not an eigenvalue of $(-\De)^s$ with homogenous Dirichlet boundary condition in $\mb R^n \setminus \Om$.}
\end{Theorem}

\begin{Theorem}\label{newthrm}
	Let $s \in (0,1)$ and $ 2s<n<4s $, then there exist $\bar{\la}>0$ such that for any $\la > \bar \la $ different from the eigenvalues of $(-\De)^s$ with homogenous Dirichlet boundary condition in $\mb R^n \setminus \Om$, \eqref{BS-10} has a nontrivial weak solution.
\end{Theorem}

\begin{Theorem}\label{thrm3}
	Let $\la <0$ and $\Om \not\equiv \mb R^n$ be a strictly star shaped bounded domain (with respect to origin) with $C^{1,1}$ boundary,  then  \eqref{BS-10} cannot have a  nonnegative nontrivial solution.
\end{Theorem}
Consider the space $X$ defined as
\[X= \left\{u|\;u:\mb R^n \to \mb R \;\text{is measurable},\;
u|_{\Om} \in L^2(\Om)\;
\text{and}\;  \frac{(u(x)- u(y))}{ |x-y|^{\frac{n}{2}+s}}\in
L^2(Q)\right\},\]
where $Q=\mb R^{2n}\setminus(\mc C\Om\times \mc C\Om)$ and
$\mc C\Om := \mb R^n\setminus\Om$ endowed with the norm
\[\|u\|_X = \|u\|_{L^2(\Om)} + \left[u\right]_X,\]
where
\[\left[u\right]_X= \left( \int_{Q}\frac{|u(x)-u(y)|^{2}}{|x-y|^{n+2s}}\,\mathrm{d}x\mathrm{d}y\right)^{\frac12}.\]
Then we define $ X_0 = \{u\in X : u = 0 \;\text{a.e. in}\; \mb R^n\setminus \Om\}$ and we have the Poincare type inequality:  there exists a constant $C>0$ such that $\|u\|_{L^{2}(\Om)} \le C [u]_X$, for all $u\in X_0$. Hence,  $\|u\|=[u]_X$ is a norm on ${X_0}$. Moreover, $X_0$ is a Hilbert space and $C_c^{\infty}(\Om)$ is dense in $X_0$. For details on these spaces and variational setup we refer to \cite{bn-serv}.
\begin{Definition}
	We say that $u \in X_0$ is a weak solution of \eqref{BS-10} if
	\begin{equation*}
	\begin{split}
	\int_{Q}\frac{(u(x)-u(y))(\varphi(x)-\varphi(y))}{|x-y|^{n+2s}}~dx dy =& \int_{\Om}\int_{\Om}\frac{|u(x)|^{2^*_{\mu,s}}|u(y)|^{2^*_{\mu,s}-2}u(y)\varphi(y)}{|x-y|^{\mu}}~dx dy \\ &+ \la \int_{\Om}u\varphi ~dx,\; \text{for every}\; \varphi \in {X_0}.
	\end{split}
	\end{equation*}
	
\end{Definition}
The corresponding energy functional associated to the problem \eqref{BS-10} is given by
\[I_{\la}(u)= I(u) := \frac{\|u\|^2}{2} - \frac{1}{22^*_{\mu,s}}\int_{\Om}\int_{\Om}\frac{|u(x)|^{2^*_{\mu,s}}|u(y)|^{2^*_{\mu,s}}}{|x-y|^{\mu}}~ dx dy
-\frac{\la}{2}\int_{\Om}|u|^2 dx.\]
Using Hardy-Littlewood-Sobolev inequality, we can show that $I \in C^1(X_0,\mb R)$ and the critical points of $I$ corresponds to weak solution {of} \eqref{BS-10}.
We define
\begin{equation*}
S^H_s := \displaystyle \inf\limits_{H^s(\mb R^n)\setminus \{0\}} \frac{\displaystyle\int_{\mb R^n}\int_{\mb R^n} \frac{|u(x)-u(y)|^2}{|x-y|^{n+2s}}{\,dxdy}}{\displaystyle\left(\int_{\mb R^n} \int_{\mb R^n}\frac{|u(x)|^{2^*_{\mu,s}}|u(y)|^{2^*_{\mu,s}}}{|x-y|^{\mu}} dx dy\right)^{\frac{1}{2^*_{\mu,s}}}}
\end{equation*}
as the best constant which is achieved if and only if $u$ is of the form
\[C\left( \frac{t}{t^2+|x-x_0|^2}\right)^{\frac{n-2s}{2}}, \; \; \text{for all}\; x \in \mb R^n,\]
for some $x_0 \in \mb R^n$, $C>0$ and $t>0$.  Moreover,
$
S^H_s ~ C(n,\mu)^{\frac{1}{2^*_{\mu,s}}}= S_s,$
where $S_s$ is the best constant of the Sobolev imbedding $H^s(\mathbb R^n)$ into $L^2(\mathbb R^n).$
Using suitable translation and dilation of the minimizing sequence, we proved-
\begin{Lemma}
	Let
	\begin{equation*}
	S^H_s(\Om) := \inf\limits_{X_0\setminus\{0\}} \frac{\displaystyle\int_{Q} \frac{|u(x)-u(y)|^2}{|x-y|^{n+2s}}{\,\mathrm{d}x\mathrm{d}y}}{\displaystyle\left(\int_{\Om} \int_{\Om}\frac{|u(x)|^{2^*_{\mu,s}}|u(y)|^{2^*_{\mu,s}}}{|x-y|^{\mu}}\mathrm{d}x\mathrm{d}y\right)^{\frac{1}{2^*_{\mu,s}}}}.
	\end{equation*}
	Then $S^H_s(\Om)= S^H_s$ and $S^H_s(\Om)$ is never achieved except when $\Om = \mb R^n$.
\end{Lemma}
Since \eqref{BS-10} has a lack of compactness due to the presence of the critical exponent, we needed a Brezis-Lieb type lemma which can be proved in the spirit of \eqref{BS-2} as follows-
\begin{equation*}
\begin{split}
\int_{\mb R^n} \int_{\mb R^n} \frac{|u_k(x)|^{2^*_{\mu,s}}|u_k(y)|^{2^*_{\mu,s}}}{|x-y|^{\mu}}~\mathrm{d}x\mathrm{d}y &- \int_{\mb R^n} \int_{\mb R^n} \frac{|(u_k-u)(x)|^{2^*_{\mu,s}}|(u_k-u)(y)|^{2^*_{\mu,s}}}{|x-y|^{\mu}}~\mathrm{d}x\mathrm{d}y\\
& \rightarrow \int_{\mb R^n} \int_{\mb R^n} \frac{|u(x)|^{2^*_{\mu,s}}|u(y)|^{2^*_{\mu,s}}}{|x-y|^{\mu}}~\mathrm{d}x \mathrm{d}y \; \text{ as}\;  k \rightarrow \infty
\end{split}
\end{equation*}
where $\{u_k\}$ is a bounded sequence in $L^{2^*_s}(\mb R^n)$ such that $u_k \rightarrow u$ almost everywhere in $\mb R^n$ as $n \rightarrow \infty$. Next, we prove the following properties concerning the compactness of Palais-Smale sequences. If $\{u_k\}$ is a Palais-Smale sequence of $I$ at $c$. Then
\begin{enumerate}
	\item[(i)] $\{u_k\}$ must be bounded in $X_0$ and its weak limit is a weak solution of \eqref{BS-10},
	\item[(ii)] $\{u_k\}$ has a convergent subsequence if
	\[c < \frac{n+2s-\mu}{2(2n-\mu)}(S^H_s)^{\frac{2n-\mu}{n+2s-\mu}}.\]
\end{enumerate}
Let us consider the sequence of eigenvalues of the operator $(-\De)^s$ with homogenous Dirichlet boundary condition in $\mb R^n \setminus \Om$, denoted by
\[0 < \la_1 < \la_2 \leq \la_3 \leq \ldots \leq \la_{j} \leq \la_{j+1}\leq \ldots\]
and $\{e_j\}_{j \in \mb N} \subset L^{\infty}(\Om)$ be the  corresponding sequence of eigen functions. We also consider this sequence of $e_j$'s to form an orthonormal basis of $L^2(\Om)$ and orthogonal basis of $X_0$. We then dealt with the cases $\la \in (0,\la_{1})$ and $\la \in (\la_{r}, \la_{r+1})$ separately. We assume $0 \in \Om$ and fix $\delta>0$ such that $B_{\delta}\subset \Om \subset B_{{\hat k}\delta}$, {for some $\hat k>1$}. Let $\eta \in C^{\infty}(\mb R^n)$ be such that $0\leq \eta \leq 1$ in $\mb R^n$, $\eta \equiv 1$ in $B_{\delta}$ and $\eta \equiv 0$ in $\mb R^n \setminus \Om$. For $\epsilon >0$, we define the function $u_\epsilon$ as follows
\[u_\epsilon(x) := \eta(x)U_{\epsilon}(x),\]
for $x \in \mb R^n$, where $ \displaystyle U_{\epsilon}(x) = \epsilon^{-\frac{(n-2s)}{2}}\; {\left(\frac{u^*\left(\frac{x}{\epsilon}\right)}{\|u^*\|_{{2^*_s}} }\right)}$ and $u^*(x)=  {\alpha\left(\beta^2 + \left|\frac{x}{ S_s^{\frac{1}{2s}} }\right|^2\right)^{-\frac{n-2s}{2}}}$ with $\alpha \in \mb R \setminus \{0\}$, $ \beta >0$.
We obtained the following important asymptotic estimates
\begin{Proposition}\label{estimates1}
	The following estimates holds true:
	\begin{equation*}\label{esti-new}
	\int_{\mb R^n}\frac{|u_\epsilon(x)- u_\epsilon(y)|^2}{|x-y|^{n+2s}}~\mathrm{d}x\mathrm{d}y \leq \left((C(n,\mu))^{\frac{n-2s}{2n-\mu}}S^H_s\right)^{\frac{n}{2s}}+ O(\epsilon^{n-2s}),
	\end{equation*}
	\begin{equation*}
	\left(\int_{\Om}\int_{\Om}\frac{|u_\epsilon(x)|^{2^*_{\mu,s}}|u_\epsilon(y)|^{2^*_{\mu,s}}}{|x-y|^{\mu}}~\mathrm{d}x\mathrm{d}y \right)^{\frac{n-2s}{2n-\mu}}\leq (C(n,\mu))^{\frac{n(n-2s)}{2s(2n-\mu)}} (S^H_s)^{\frac{n-2s}{2s}}+ {O(\epsilon^{n})},
	\end{equation*}
	and
	\begin{equation*}
	\left(\int_{\Om}\int_{\Om}\frac{|u_\epsilon(x)|^{2^*_{\mu,s}}|u_\epsilon(y)|^{2^*_{\mu,s}}}{|x-y|^{\mu}}~\mathrm{d}x\mathrm{d}y \right)^{\frac{n-2s}{2n-\mu}}\geq {\left((C(n,\mu))^{\frac{n}{2s}} (S^H_s)^{\frac{2n-\mu}{2s}}- O\left(\epsilon^{n}\right)\right)^{\frac{n-2s}{2n-\mu}}.}
	\end{equation*}
\end{Proposition}
When $n \geq 4s$, we proved that the energy functional $I_\la$ satisfies the Mountain pass geometry if $\la \in (0,\la_1)$ and Linking Theorem geometry if $\la \in (\la_r,\la_{r+1})$. Also in both the cases, Proposition \ref{estimates1} helped us to show that for small enough $\e>0$
\begin{equation}\label{elem1}
\frac{\|u_\e\|^2 - \la \int_{\Om}{|u_\e|^2}\mathrm{d}x}{\left(\int_{\Om}\int_{\Om}\frac{|u_\e(x)|^{2^*_{\mu,s}}{|u_\e(y)|^{2^*_{\mu,s}}}}{|x-y|^{\mu}}~\mathrm{d}x\mathrm{d}y \right)^{\frac{n-2s}{2n-\mu}}} < S^H_s.
\end{equation}
Then the proof of Theorem \ref{thrm1} follows by applying Mountain Pass Lemma and Linking Theorem. On the other hand when $2s<n<4s$, \eqref{elem1} could be proved only when $\la > \bar \la$ for some suitable $\bar \la>0$, when $\e>0$ is sufficiently small. Hence again applying Mountain Pass Lemma and Linking Theorem in this case too, we prove Theorem \ref{newthrm}.
To prove Theorem \ref{thrm3}, we first prove that if $\la <0$ then any solution $u \in X_0$ of \eqref{BS-10} belongs to $L^\infty(\Om)$ which implied that when $\Om$ is a $C^{1,1}$ domain then $u/\delta^s \in C^\alpha(\bar \Om)$
for some $\alpha >0$ (depending on $\Om$ and $s$) satisfying $\alpha < \min\{s,1-s\}$, where $\delta(x) = \text{dist}(x, \partial \Om)$ for $x\in \Om$. Then using $(x.\nabla u)$ as a test function in \eqref{BS-10}, we proved the following Pohozaev type identity-
\begin{Proposition}\label{Poho}
	If  $\la <0$, $\Om$ be bounded $C^{1,1}$ domain and $u \in L^{\infty}(\Om)$ solves \eqref{BS-10}, then
	\begin{equation*}
	\begin{split}
	\frac{2s-n}{2}&\int_{\Om}u(-\De)^su~\mathrm{d}x - \frac{\Gamma(1+s)^2}{2}\int_{\partial \Om}\left(\frac{u}{\delta^s}\right)^2(x.\nu)\mathrm{d}\sigma\\
	&={-}\left(\frac{2n-\mu}{22^*_{\mu,s}}\int_{\Om}\int_{\Om}\frac{|u(x)|^{2^*_{\mu,s}}|u(y)|^{2^*_{\mu,s}}}{|x-y|^{\mu}}~\mathrm{d}x\mathrm{d}y+ {\frac{\la n}{2}} \int_{\Om}|u|^2\mathrm{d}x\right),
	\end{split}
	\end{equation*}
	where $\nu$ denotes the unit outward normal to $\partial \Om$ at $x$ and $\Gamma$ is the Gamma function.
\end{Proposition}
Using Proposition \ref{Poho}, Theorem \ref{thrm3} easily followed.
\subsection{Magnetic Choquard equations}
Very recently L\"{u} \cite{Lu} studied the problem
\begin{equation}\label{intro4}
(-i \nabla +A(x))^2u + (g_0+\mu g)(x) u = (|x|^{-\alpha}* |u|^p)|u|^{p-2}u, \; u \in H^1(\mb R^n, \mb C),
\end{equation}
where $n\geq 3$, $\alpha \in (0,n)$, $p \in \left( \frac{2n-\alpha}{n}, \frac{2n-\alpha}{n-2}\right)$, $A = (A_1, A_2, \ldots, A_n): \mb R^n \rightarrow {\mb R^n}$ is a vector (or magnetic) potential such that {$A \in L^n_{\text{loc}}(\mb R^n, \mb R^n)$} and {$A$ is continuous at ${0}$},  $g_0$ and $g$ are real valued functions on $\mb R^n$ satisfying some {necessary} conditions and $\mu >0$.  He proved the existence of ground state solution {when} $\mu \geq \mu^*$, for some $\mu^*>0$ and concentration behavior of solutions as $\mu \rightarrow \infty$. Salazar {in} \cite{szr} showed existence of vortex type solutions for the stationary nonlinear magnetic Choquard equation
\begin{equation*}
(-i \nabla +A(x))^2u + W(x)u = (|x|^{-\alpha}* |u|^p)|u|^{p-2}u \; \text{in} \; \mb R^n,
\end{equation*}
where  $p \in \left[2, 2^*_\alpha \right)$ and $W: \mb R^n \to  \mb R$ is bounded electric potential. Under some assumptions on decay of $A$ and $W$ at infinity, Cingloni, Sechi and Squassina in \cite{css} showed existence of family of solutions.  Schr\"{o}dinger equations with magnetic field and  Choquard type nonlinearity has also been studied in \cite{MSS,MSS1}. But the critical case in \eqref{intro4} was still open which motivated us to study the problem $(P_{\la,\mu})$ in \cite{TS-2}:
$$(P_{\la,\mu})
\left\{
\begin{array}{ll}
(-i \nabla+A(x))^2u + \mu g(x)u = \la u + (|x|^{-\alpha} * |u|^{2^*_\alpha})|u|^{2^*_\alpha-2}u & \mbox{in } \mb R^n \\
u \in H^1(\mb R^n, \mb C)
\end{array}
\right.$$
where $n \geq 4, 2^{*}_{\alpha}= \frac{2n-\alpha}{n-2}, \alpha \in (0,n), \mu>0, \la>0, A = (A_1, A_2, \ldots, A_n): \mathbb{R}^n \rightarrow \mathbb{R}^n$ is a vector(or magnetic) potential such that
$A\in L^{n}_{loc}(\mathbb{R}^n, \mathbb{R}^n)$ and $A$ is continuous at ${0}$ and $g(x)$ satisfies the following assumptions:
\begin{enumerate}
	\item[(g1)] $g \in C(\mb R^n,\mb R)$, $g \geq 0$ and $\Om := \text{interior of}\;g^{-1}(0)$ is  a nonempty bounded set with smooth boundary and $\overline{\Om}= g^{-1}(0)$.
	\item[(g2)] There exists $M>0$ such that ${\mc L}\{x \in \mb R^n:\; g(x)\leq M\} < {+\infty}$, where ${\mc L}$ denotes the Lebesgue measure in $\mb R^n$.
\end{enumerate}
Let us define $-\nabla_A := (-i\nabla +A)$ and
$$ H^1_A(\mb R^n, \mb C)= \left\{u \in L^2(\mb R^n,\mb C) \;: \; \nabla_A u \in L^2(\mb R^n, \mb C^n)\right\}.$$
Then $H^1_A(\mb R^n, \mb C)$ is a Hilbert space with the inner product
$$\langle u,v\rangle_A = \text{Re} \left(\int_{\mb R^n} (\nabla_A u \overline{\nabla_A v} + u \overline v)~\mathrm{d}x \right),$$
where $\text{Re}(w)$ denotes the real part of $w \in \mb C$ and $\bar w$ denotes its complex conjugate. The associated norm $\|\cdot\|_A$ on the space $H^1_A(\mb R^n, \mb C)$ is given by
$$\|u\|_A= \left(\int_{\mb R^n}(|\nabla_A u|^2+|u|^2)~\mathrm{d}x\right)^{\frac{1}{2}}.$$
We call $H^1_A(\mb R^n, \mb C)$ simply $H^1_A(\mb R^n)$. Let $H^{0,1}_A(\Om, \mb C)$ (denoted by $H^{0,1}_A(\Om)$ for simplicity) be the Hilbert space defined by the closure of $C_c^{\infty}(\Om, \mb C)$ under the scalar product {$\langle u,v \rangle_A= \textstyle\text{Re}\left(\int_{\Om}(\nabla_A u \overline{\nabla_A v}+u \overline v)~\mathrm{d}x\right)$}, where $\Om = \text{interior of } g^{-1}(0)$. Thus {norm on $H^{0,1}_A(\Om)$ is given by}
\[\|u\|_{H^{0,1}_A(\Om)}= \left(\int_\Om (|\nabla_A u|^2+|u|^2)~\mathrm{d}x\right)^{\frac{1}{2}}.\]
Let $E = \left\{u \in H^1_A(\mb R^n): \int_{\mb R^n}g(x)|u|^2~\mathrm{d}x < +\infty\right\}$ be the Hilbert space equipped with the inner product
$$\langle u,v\rangle = \text{Re} \left(\int_{\mb R^n}\left(\nabla_A u \overline{\nabla_A v}~\mathrm{d}x + g(x)u\bar v \right)~\mathrm{d}x\right)$$
and the associated norm
$\|u\|_E^2=  \int_{\mb R^n}\left(|\nabla_A u|^2+ g(x)|u|^2\right)~\mathrm{d}x. $
Then $\|\cdot\|_E$ is clearly equivalent to each of the norm
$\|u\|_\mu^2= \int_{\mb R^n}\left(|\nabla_A u|^2+ \mu g(x)|u|^2\right)~\mathrm{d}x$
for $\mu >0$. We have the following well known \textit{diamagnetic inequality} (for detailed proof, see \cite{LL}, Theorem $7.21$ ).
\begin{Theorem}\label{dia_eq}
	If $u \in H^1_A(\mb R^n)$, then $|u| \in H^1(\mb R^n,\mb R)$ and
	$$|\nabla |u|(x)| \leq |\nabla u(x)+ i A(x)u(x)| \; \text{for a.e.}\; x \in \mb R^n.$$
\end{Theorem}
So for each $q \in [2,2^*]$, there exists constant $b_q>0$ (independent of $\mu$) such that
\begin{equation}\label{mg-eq1}
|u|_q \leq b_q\|u\|_\mu, \; \text{for any}\; u \in E,
\end{equation}
where $|\cdot|_q$ denotes the norm in $L^q(\mb R^n,\mb C)$ and $2^*= \textstyle\frac{2n}{n-2}$ is the Sobolev critical exponent. {Also $H^1_A(\Om) \hookrightarrow L^q(\Om, \mb C)$ is continuous for each $1\leq q \leq 2^*$ and compact when $1\leq q < 2^*$.} We denote ${\la_1(\Om)}>0$ as the best constant of the  embedding ${H^{0,1}_A(\Om)}$ into $L^2(\Om, \mb C) $ given by
\[ \la_1(\Om) = \inf\limits_{u \in H^{0,1}_A(\Om)}\left\{\int_{\Om}|\nabla_A u|^2 ~\mathrm{d}x : \; {\int_{\Om}|u|^2~\mathrm{d}x}=1\right\}\]
which is also the first eigenvalue of $-\Delta_A := (-i\nabla +A)^2$ on $\Om$ with boundary condition $u=0$. In \cite{TS-2}, we consider the problem
$$ (P_\la)\;\;
(-i\nabla +A(x))^2u = \la u + (|x|^{-\alpha}*|u|^{2^*_\alpha})|u|^{2^*_\alpha-2}u \; \mbox{in } \;\Om, \;\;
u=0 \; \mbox{on }\; \partial \Om
$$
and proved the following main results:
\begin{Theorem}\label{MT1}
	For every $\la \in (0, \la_1(\Om))$ there exists a $\mu(\la)>0$ such that $(P_{\la,\mu})$ has a least energy solution $u_\mu$ for each $\mu\geq \mu(\la)$.
\end{Theorem}

\begin{Theorem}\label{MT2}
	Let $\{u_m\}$ be a sequence of non-trivial solutions of $(P_{\la,\mu_m})$ with $\mu_m \rightarrow \infty$ and $I_{\la,\mu_m}(u_m) \rightarrow c< \frac{n+2-\alpha}{2(2n-\alpha)}S_A^{\frac{2n-\alpha}{n+2-\alpha}}$ as $m \rightarrow \infty$. Then $u_{m}$ concentrates at a solution of $(P_\la)$.
\end{Theorem}
We give some definitions below-
\begin{Definition}
	We say that a function $u \in {E}$ is a weak solution of $(P_{\la,\mu})$ if
	\[\text{Re}\left( \int_{\mb R^n} \nabla_A u \overline{\nabla_A v}~\mathrm{d}x + \int_{\mb R^n}(\mu g(x)-\la )u\overline{v}~\mathrm{d}x- \int_{\mb R^n} (|x|^{-\alpha}* |u|^{2^*_\alpha})|u|^{2^*_\alpha-2}u \overline{v}~\mathrm{d}x \right)=0\]
	for all $v \in {E}$.
\end{Definition}
\begin{Definition}
	A solution $u$ {of $(P_{\la,\mu})$} is said to be a least energy solution if the energy functional
	\[I_{\la,\mu}(u)=  \int_{\mb R^n}\left( \frac12\left(|\nabla_A u|^2+ (\mu g(x)-\la)|u|^2\right) - \frac{1}{22^*_\alpha}(|x|^{-\alpha}* |u|^{2^*_\alpha})|u|^{2^*_\alpha}\right) ~\mathrm{d}x\]
	achieves its minimum at $u$ over all the nontrivial solutions of $(P_{\la,\mu})$.
\end{Definition}
\begin{Definition}
	A sequence of solutions $\{u_k\}$ of $(P_{\la,\mu_k})$ is said to concentrate at a solution $u$ of $(P_\la)$ if a subsequence converges strongly to $u$ in $H^1_A(\mb R^n)$ as $\mu_k \rightarrow \infty$.
\end{Definition}
We first proved the following Lemma.
\begin{Lemma}\label{comp_lem1}
	Suppose $\mu_m \geq 1$ and $u_m \in E$ be such that $\mu_m \rightarrow \infty$ as $m \rightarrow \infty$ and {there exists a $K>0$ such that} $\|u_m\|_{\mu_m} < K$, for all $m \in \mb N$. Then there exists a $u \in H^{0,1}_A(\Om)$ such that (upto a subsequence), $u_m \rightharpoonup u$ weakly in $E$ and $u_m \rightarrow u$ strongly in $L^2(\mb R^n)$  as $m \rightarrow \infty$.
\end{Lemma}
Then we define an operator $T_\mu := -\Delta_A + \mu g(x)$ on $E$ given by
\[\big(T_\mu(u),v\big)= \text{Re}\left(\int_{\mb R^n}(\nabla_A u \overline{\nabla_A v}+ \mu g(x)u\overline v)~\mathrm{d}x\right).\]
Clearly $T_\mu$ is a self adjoint operator and if $a_\mu := \inf \sigma(T_\mu)$, i.e. the infimum of the spectrum of $T_\mu$, then $a_\mu$ can be characterized as
\[0 \leq a_\mu = \inf \{\big(T_\mu(u),u\big): \; u \in E,\; \|u\|_{L^2}=1\}= \inf \{\|u\|_\mu^2:\; u\in E,\;\|u\|_{L^2}=1\}.\]
Then considering a minimizing sequence of $a_\mu$, we were able to prove that
for each $\la \in (0, \la_1(\Om))$, there exists a $\mu(\la)>0$  such that $a_\mu \geq (\la +\la_1(\Om))/2$ whenever $\mu \geq \mu(\la)$. As a consequence
\[\big((T_\mu-\la)u,u) \geq \beta_\la \|u\|_\mu^2\]
for all $u \in E$, $\mu \geq \mu(\la)$, where $\beta_\la := (\la_1(\Om)-\la)/(\la_1(\Om)+\la)$. We fix $\la \in (0, \la_1(\Om))$ and $\mu \geq \mu(\la)$. Using standard techniques, we established the following concerning any Palais Smale sequence $\{u_k\}$ of $I_{\la,\mu}$ -
\begin{enumerate}
	\item[(i)] $\{u_m\}$ must be bounded in $E$ and its weak limit is a solution of $(P_{\la,\mu})$,
	\item[(ii)] $\{u_m\}$ has a convergent subsequence when $c$ satisfies
	\[c \in \left(-\infty,  \frac{n+2-\alpha}{2(2n-\alpha)} S_A^{\frac{2n-\alpha}{n+2-\alpha}}\right)\]
	where $S_A$ is defined as follows
	\[S_A = {\inf_{u \in H^1_A(\mb R^n) \setminus \{0\}} \frac{\displaystyle \int_{\mb R^n}|\nabla_Au|^2~\mathrm{d}x}{\displaystyle \int_{\mb R^n}(|x|^{-\alpha}* |u|^{2^*_\alpha})|u|^{2^*_\alpha}~\mathrm{d}x}}.\]
\end{enumerate}
Using asymptotic estimates using the family
$U_\epsilon (x)= (n(n-2))^{\frac{n-2}{4}}\left(\frac{\epsilon}{\epsilon^2+|x|^2}\right)^{\frac{n-2}{4}}$,
we showed that-
\begin{Theorem}\label{S_Aattain}
	If $g\geq 0$ and $A \in L^n_{\text{loc}}(\mb R^n, \mb R^n)$, then the infimum $S_A$ is attained if and only if $\text{curl }A \equiv 0$.
\end{Theorem}
Our next step was to introduce the Nehari manifold
\[\mc N_{\la,\mu}= \left\{u \in E\setminus \{0\}:\; \langle I^\prime_{\la,\mu}(u),u\rangle=0\right\}\]
and consider the minimization problem
$k_{\la,\mu}:= \inf_{u \in \mc N_{\la,\mu}} I_{\la,\mu}(u)$.
Using the family $\{U_\e\}$, we showed that
\[k_{\la,\mu}  < \frac{n+2-\alpha}{2(2n-\alpha)} S_A^{\frac{2n-\alpha}{n+2-\alpha}}.\]
Then the proof of Theorem \ref{MT1} followed by using the Ekeland Variational Principle over $\mc N_{\la,\mu}$. The proof of Theorem \ref{MT2} followed from Lemma \ref{comp_lem1} and the Brezis-Lieb type lemma for the Riesz potentials.
\begin{Remark}
	These  results can be generalized to the problems involving  fractional magnetic operators:
	\begin{equation*}
	(P_{\la,\mu}^s)\left\{
	\begin{array}{rlll}
	& (-\De)^s_A u + \mu g(x)u = \la u + (|x|^{-\alpha} * |u|^{2^*_{\alpha,s}})|u|^{2^*_{\alpha,s}-2}u  \;\text{in} \; \mb R^n ,\\
	& u \in H^s_A(\mb R^n, \mb C)
	\end{array}
	\right.
	\end{equation*}
	where $n \geq 4s$, $s \in (0,1)$ and $\alpha\in (0,n)$. Here $2^*_{\alpha,s}=\textstyle \frac{2n-\alpha}{n-2s}$ is the critical exponent in the sense of Hardy-Littlewood-Sobolev inequality. We assume the same conditions on $A$ and $g$ as before. For $u \in C_c^\infty(\Om)$, the fractional magnetic operator $(-\De)^s_A$, up to a normalization constant, is defined by
	\[(-\De)^s_A u (x) = 2\lim_{\e \to 0^+} \int_{\mb R^n \setminus B_{\e}(x)} \frac{u(x)-e^{i(x-y)\cdot A\left(\frac{x+y}{2}\right)}u(y) }{|x-y|^{n+2s}}\mathrm{d}y \]
	for all $x \in \mb R^n$. With proper functional setting, we can prove the existence and concentration  results for the problem $(P_{\la,\mu}^s)$ employing the same arguments as in the local magnetic operator case.
\end{Remark}

\subsection{Singular problems involving Choquard nonlinearity}
The paper by Crandal, Rabinowitz and Tartar \cite{crt} is the starting point on semilinear problem with singular nonlinearity. A lot of work has been done related to the existence and multiplicity results for singular problems, see \cite{haitao, hirano2, hirano1}. Using splitting Nehari manifold technique, authors in \cite{hirano1} studied the existence of multiple solutions of the problem:
\begin{equation}\label{hr}
-\De u = \la u^{-q}+ u^p,\; u>0 \; \text{in}\;
\Om, \quad u = 0 \; \mbox{on}\; \partial\Om,
\end{equation}
where $\Om$ is smooth bounded domain in $\mb R^n$, $n\geq 1$, $p=2^*-1$, $\la>0$  and $0<q<1$. In \cite{haitao}, Haitao studied the equation \eqref{hr} for $n\geq 3$, $1<p\leq 2^{*}-1$ and showed the  existence  of two positive  solutions for maximal interval of the parameter $\la$ using monotone iterations and mountain pass lemma. But the singular problem involving Choquard nonlinearity was completely open until we studied the following problem in \cite{TS-3}
\begin{equation*}
(P_{\la}): \quad
\quad -\De u = \la u^{-q} + \left( \int_{\Om}\frac{|u(y)|^{2^*_{\mu}}}{|x-y|^{\mu}}\mathrm{d}y \right)|u|^{2^*_{\mu}-2}u, \; u>0 \; \text{in}\;
\Om, \; u = 0 \; \mbox{on}\; \partial\Om,
\quad
\end{equation*}
where $\Om \subset \mb R^n$, $n>2$ be a bounded domain with smooth boundary $\partial \Om$, $\la >0,\; 0 < q < 1, $ $ 0<\mu<n$ and  $2^*_\mu=\frac{2n-\mu}{n-2}$. The main difficulty in treating $(P_\la)$ is the presence of singular nonlinearity along with critical exponent in the sense of Hardy-Littlewood-Sobolev inequality which is nonlocal in nature. The energy functional no longer remains differentiable due to presence of singular nonlinearity, so usual minimax theorems are not applicable. Also the critical exponent term being nonlocal adds on the difficulty to study the Palais-Smale level around a nontrivial critical point.

\begin{Definition}
	We say that  $u\in H^1_0(\Om)$ is a positive weak solution of $(P_\la)$
	if $u>0$ in $\Om$ and
	\begin{equation}\label{BS-11}
	\int_{\Om} (\nabla u \nabla \psi -\la u^{-q}\psi)~\mathrm{d}x  - \int_{\Om}\int_{\Om}\frac{|u(x)|^{2^*_{\mu}}|u(y)|^{2^*_{\mu}-2}u(y)\psi(y)}{|x-y|^{\mu}}~\mathrm{d}x\mathrm{d}y = 0
	\end{equation}
	{for all} $\psi \in C^{\infty}_c(\Om)$.
\end{Definition}
We define the functional associated to $(P_\la)$ as $I : H^1_{0}(\Om) \rightarrow (-\infty, \infty]$ by
\[ I(u) = \frac12 \int_{\Om}|\nabla u|^2~ \mathrm{d}x- \frac{\la}{1-q} \int_\Om |u|^{1-q} \mathrm{d}x - \frac{1}{22^*_{\mu}}\int_{\Om}\int_{\Om}\frac{|u(x)|^{2^*_{\mu}}|u(y)|^{2^*_{\mu}}}{|x-y|^{\mu}}~\mathrm{d}x\mathrm{d}y, \]
for $u \in H^1_{0}(\Om)$. For each $0 < q <1$, we set
$ H_+ = \{ u \in H^1_0(\Om) : u \geq 0\}$ and
$$H_{+,q} = \{ u \in H_+ : u \not\equiv 0,\; |u|^{1-q} \in L^1(\Om)\} = H_+ \setminus \{0\} .$$
For each $u\in H_{+,q}$ we define the fiber map $\phi_u:\mb R^+ \rightarrow \mb R$ by $\phi_u(t)=I_\la(tu)$. Then we proved the following:
\begin{Theorem} \label{mainthrm1}
	Assume $0<q < 1$ and let $\Lambda$ be a constant defined by
	\begin{equation*}
	\begin{split}
	\Lambda = & \sup \left\{\la >0:\text{ for each}  \; u\in H_{+,q}\backslash\{0\}, ~\phi_u(t)~  \text{has two critical points in}  ~(0, \infty)\right.\\
	&\left.\text{and}\; \sup\left\{ \int_{\Om}|\nabla u|^2~\mathrm{d}x\; : \; u \in H_{+,q}, \phi^{\prime}_u(1)=0,\;\phi^{\prime\prime}_u(1)>0 \right\} \leq (2^*_\mu S_{H,L}^{2^*_\mu})^{\frac{1}{2^*_\mu-1}} \right\}.
	\end{split}
	\end{equation*}
	Then $\La >0$.
\end{Theorem}
Using the variational methods on the Nehari manifold, we proved the following multiplicity result.
\begin{Theorem}\label{mainthrm2}
	For all $\la \in (0, \Lambda)$, $(P_\la)$ has two positive weak solutions $u_\la$ and $v_\la$ in $C^\infty(\Om)\cap L^{\infty}(\Om)$.
\end{Theorem}
We also have that if $u$ is a positive weak solution of $(P_\la)$, then $u$ is a classical solution in the sense that $u \in C^\infty(\Om) \cap C(\bar \Om)$.
We define $\delta : \Om \rightarrow [0,\infty)$ by $\delta(x)=\inf\{|x-y|: y \in \partial \Om\}$, for each $x \in \Om$.
\begin{Theorem}\label{mainthrm3}
	Let $u$ be a positive weak solution of $(P_\la)$, then there exist $K,\;L>0$ such that $L\delta \leq u \leq K\delta$ in $\Om$.
\end{Theorem}
We define the Nehari manifold
\[ \mc N_{\la} = \{ u \in H_{+,q} | \left\langle I^{\prime}(u),u\right\rangle = 0 \}\]
and show that
$I$ is coercive and bounded below on $\mc N_{\la}$. It is easy to see that the points in $\mc N_{\la}$ are corresponding to critical points of $\phi_{u}$ at $t=1$. So, we divided $\mc N_{\la}$ in three sets corresponding to local minima, local maxima and points of inflexion
\begin{align*}
\mc N_{\la}^{+} &=  \{ t_0u \in \mc N_{\la} |\; t_0 > 0,~ \phi^{\prime}_u (t_0) = 0,~ \phi^{\prime \prime}_u(t_0) > 0\},\\
\mc N_{\la}^{-} = & \{ t_0u \in \mc N_{\la} |\; t_0 > 0, ~ \phi^{\prime}_u (t_0) = 0, ~\phi^{\prime \prime}_u(t_0) < 0\}
\end{align*}
and $ \mc N_{\lambda}^{0}= \{ u \in \mc N_{\la} | \phi^{\prime}_{u}(1)=0,\; \phi^{\prime \prime}_{u}(1)=0 \}$. We aimed at showing that the minimizers of $I$ over $\mc N^+$ and $\mc N^-$ forms a weak solution of $(P_\la)$. We briefly describe the key steps to show this. Using the fibering map analysis, we proved that there exist $\la_*>0$ such that for each $u\in H_{+,q}\backslash\{0\}$, there is unique $t_1$ and $t_2$ with the property that $t_1<t_2$, $t_1 u\in \mc N_{\la}^{+}$ and $ t_2 u\in \mc N_{\la}^{-}$, for all $\la \in (0,\la_*)$. This implied Theorem \ref{mainthrm1}. Also $\mc N_{\la}^{0} = \{0\}$ for all $ \la \in (0, \la_*)$. Then it is easy to see that $\sup \{ \|u\|: u \in \mc N_{\la}^{+}\} < \infty $
and $\inf \{ \|v\|: v \in \mc N_{\la}^{-} \} >0$. Suppose $u$ and $v$ are minimizers of $I$ on $\mc N_{\la}^{+}$ and $\mc N_{\la}^{-}$ respectively. Then
for each $w \in H_{+}$, we showed $u^{-q}w, v^{-q} w \in L^{1}(\Om)$ and
\begin{align}
&\int_{\Om} (\nabla u \nabla w-\la u^{-q}w)~\mathrm{d}x  - \int_{\Om}\int_{\Om}\frac{|u(y)|^{2^*_{\mu}}|u(x)|^{2^*_{\mu}-2}u(x)w(x)}{|x-y|^{\mu}}~\mathrm{d}y\mathrm{d}x \geq 0 , \label{upos}\\
&\int_{\Om} (\nabla v \nabla w -\la v^{-q}w)~\mathrm{d}x  - \int_{\Om}\int_{\Om}\frac{|v(y)|^{2^*_{\mu}}|u(x)|^{2^*_{\mu}-2}v(x)w(x)}{|x-y|^{\mu}}~\mathrm{d}y\mathrm{d}x \geq 0.\label{vpos}\end{align}
Particularly, $u,\; v >0$ almost everywhere in $\Om$. Then the claim followed using the Gat\'{e}aux differentiability of $I$. Lastly, the proof of Theorem \ref{mainthrm2} followed by proving that $I$ achieves its minimum over the sets $\mc N^+_\la$ and $\mc N^-_\la$.\\
In the regularity section, firstly, we showed that \eqref{BS-11} holds for all $\psi \in H_0^1(\Om)$ and each positive weak solution $u$ of $(P_\la)$ belongs to $L^\infty(\Om)$. Under the assumption that there exist $a\geq 0$, $R \geq 0$ and $q \leq s <1$ such that
$\De \delta \leq R\delta^{-s} \; \text{in} \; \Om_a:=\{x\in \Om, \de(x)\le a\},$
using appropriate test functions, we proved that
there exist $K>0$ such that $u \leq K\delta$ in $\Om$. To get the lower bound on $u$ with respect to $\delta$, following result from \cite{brenir} plays a crucial role.
\begin{Lemma}
	Let $\Om$ be a bounded domain in $\mb R^n$ with smooth boundary $\partial \Om$. Let $u \in L^1_{\text{loc}}(\Om)$ and assume that for some $k \geq 0$, $u$ satisfies, in the sense of distributions
	\[
	-\De u + ku \geq 0 \; \text{in} \; \Om,\quad
	u \geq 0 \;  \text{in}\; \Om.
	\]
	Then either $u \equiv 0$, or there exists $\e>0$ such that
	$u(x) \geq \e \delta(x), \; x \in \Om.$
\end{Lemma}
Additionally, we also prove that the solution can be more regular in a restricted range of $q$.
\begin{Lemma}
	Let $q\in (0,\frac{1}{n})$ and  let $u \in H^1_0(\Om)$ be a positive weak solution of $(P_\la)$, then $u \in C^{1+\alpha}(\bar \Om)$ for some $0<\alpha<1$.
\end{Lemma}

\section{System of equations with Choquard type nonlinearity}
In this section, we briefly illustrate some existence and multiplicity results proved concerning system of Choquard equations with nonhomogeneous terms. We consider the nonlocal operator that is the fractional Laplacian and since the Choquard nonlinearity is also a nonlocal one, such problems are often called 'doubly nonlocal problems'.  We employ the method of Nehari manifold to achieve our goal.
\subsection{Doubly nonlocal $p$-fractional Coupled elliptic system}
The $p$-fractional Laplace operator is defined  as
\[(-\De)^s_pu(x)= 2 \lim_{\e \searrow 0} \int_{|x|>\e} \frac{|u(x)-u(y)|^{p-2}(u(x)-u(y))}{|x-y|^{n+sp}}~dy, \; \forall x\in \mb R^n,\]
which is nonlinear and nonlocal in nature. This definition matches to linear fractional Laplacian operator $(-\De)^s$, {up to} a normalizing constant depending on $n$ and $s$, when $p=2$. {The operator} $(-\De)^s_p$ is degenerate when $p>2$ and singular when $1<p<2$. For details, refer \cite{hitch}. Our concern lies in the nonhomogenous Choquard equations and system of equations. Recently,  authors in \cite{XXW} and \cite{yang-zamp} showed multiplicity of positive solutions for  a nonhomogeneous Choquard equation using Nehari manifold. The motivation behind such problems lies in the famous article by  Tarantello \cite{tarantello} where author used the structure of associated Nehari manifold to obtain the multiplicity of solutions  for the following nonhomogeneous Dirichlet problem on bounded domain $\Om$
\begin{equation*}
-\De u = |u|^{2^*-2}u+f \;\text{in}\; \Om,\; u=0 \;\text{on}\; \partial \Om.
\end{equation*}
System of elliptic equations involving  $p$-fractional Laplacian {has} been studied in \cite{CD,CS}  using Nehari manifold techniques. Very recently, Guo et al. \cite{guo} studied  a nonlocal system involving fractional Sobolev critical exponent and fractional Laplacian. There are not many results on elliptic systems with non-homogeneous nonlinearities in the literature but we cite \cite{choi,faria,ww} as some very recent works on the study of fractional elliptic systems.

\noi Motivated by these articles, we consider the following nonhomogenous quasilinear system of equations with perturbations involving $p$-fractional Laplacian in \cite{TS-4}:\\
\noi Let $p\geq 2, s\in (0,1), n>sp$, $\mu \in (0,n)$, $\frac{p}{2}\left( 2-\frac{\mu}{n}\right) < q <\frac{p^*_s}{2}\left( 2-\frac{\mu}{n}\right)$, $\alpha,\beta,\gamma >0$,
\begin{equation*}
(P)\left\{
\begin{array}{rlll}
(-\De)^s_p u+ a_1(x)u|u|^{p-2} &= \alpha(|x|^{-\mu}*|u|^q)|u|^{q-2}u+ \beta (|x|^{-\mu}*|v|^q)|u|^{q-2}u\\
& \quad \quad + f_1(x)\; \text{in}\; \mb R^n,\\
(-\De)^s_p v+ a_2(x)v|v|^{p-2} &= \gamma(|x|^{-\mu}*|v|^q)|v|^{q-2}v+ \beta (|x|^{-\mu}*|u|^q)|v|^{q-2}v\\
&\quad \quad + f_2(x)\; \text{in}\; \mb R^n,
\end{array}
\right.
\end{equation*}
where  $0< a_i \in C^1(\mb R^n, \mb R)$, $i=1,2$ and $f_1,f_2: \mb R^n \to \mb R$ are perturbations. Here
$p^*_s = \frac{np}{n-sp}$
is the critical exponent associated with the embedding of the fractional Sobolev space $W^{s,p}(\mb R^n)$ into $L^{p_s^*}(\mb R^n).$
Wang et. al in {\cite{wangetal}} studied the problem $(P)$ in the local case $s=1$ and obtained a partial multiplicity results. We improved their results and showed the multiplicity results with a weaker assumption \eqref{star0} of $f_1$ and $f_2$ below. For $i=1,2$ we introduce the spaces
\[Y_i:= \left\{u \in W^{s,p}(\mb R^n): \; \int_{\mb R^n}a_i(x)|u|^p~dx < +\infty \right\}\]
then $Y_i$ are Banach spaces equipped with the norm
\[\|u\|_{Y_i}^p = \int_{\mb R^n}\int_{\mb R^n}\frac{|u(x)-u(y)|^p}{|x-y|^{n+sp}}dxdy+ \int_{\mb R^n}a_i(x)|u|^pdx.  \]
We define the product space
$Y= Y_1 \times Y_2$ which { is a reflexive Banach space} with the norm
\[\|(u,v)\|^p := \|u\|_{Y_1}^p+ \|v\|_{Y_2}^p, \]
for all $(u,v)\in Y$.
We  assume the following condition on $a_i$, for $i=1,2$
\begin{equation*}\label{cond-on-lambda}
(A)\;\; a_i \in C(\mb R^n),\; a_i >0 \; \text{and there exists}\; M_i>0 \; \text{such that}\; \mu(\{x\in \mb R^n: a_i \leq M_i\})< \infty.
\end{equation*}
Then under the condition (A) on $a_i$, for $i=1,2$, we get $Y_i$ is continuously imbedded into $ L^r(\mb R^n)$  for $r \in [p,p^*_s]$.
To obtain our results, we assumed the following condition on perturbation terms:
\begin{equation}\label{star0}
\int_{\mb R^n} (f_1u+ f_2 v)< C_{p,q}\left(\frac{2q+p-1}{4pq}\right)\|(u,v)\|^{\frac{p(2q-1)}{2q-p}}
\end{equation}
for all $(u,v)\in Y$ such that $$\int_{\mb R^n}\left(\alpha(|x|^{-\mu}*|u|^q)|u|^q +2\beta (|x|^{-\mu}*|u|^q)|v|^q
+\gamma (|x|^{-\mu}*|v|^q)|v|^q \right)dx= 1$$ and
\[C_{p,q}= \left(\frac{p-1}{2q-1}\right)^{\frac{2q-1}{2q-p}}\left(\frac{2q-p}{p-1}\right).\]
It is easy to see that
$2q > p\left(\frac{2n-\mu}{n}\right) > p-1> \frac{p-1}{2p-1}$
which implies
$\frac{2q+p-1}{4pq}<1.$
So \eqref{star0} implies that
\begin{equation}\label{star00}
\int_{\mb R^n} (f_1u+ f_2 v)< C_{p,q}\|(u,v)\|^{\frac{p(2q-1)}{2q-p}}
\end{equation}
which we used more frequently rather than our actual assumption \eqref{star0}.
Now, the main results goes as follows.
\begin{Theorem}\label{mainthrm}
	Suppose
	$\displaystyle\frac{p}{2}\left(\frac{2n-\mu}{n}\right)< q < \displaystyle\frac{p}{2}\left(\frac{2n-\mu}{n-sp}\right)$,
	$\mu \in (0,n)$ and $(A)$ holds true. Let $0 \not \equiv f_1,f_2 \in L^{\frac{p}{p-1}}(\mb R^n)$ satisfies \eqref{star0}
	then $(P)$ has at least two weak solutions, in  which one forms a local minimum of $J$ on $Y$.
	Moreover if $f_1,f_2 \geq 0$ then this solution is a nonnegative weak solution.
\end{Theorem}
If $u,\phi \in W^{s,p}(\mb R^n)$, we use the notation $\langle u,\phi\rangle $ to denote
\[\langle u,\phi\rangle := \int_{\mb R^n}\int_{\mb R^n} \frac{(u(x)-u(y))|u(x)-u(y)|^{p-2}(\phi(x)-\phi(y))}{|x-y|^{n+sp}}dxdy. \]
\begin{Definition}
	A pair of functions $(u,v)\in Y$ is said to be a weak solution to $(P)$ if
	\begin{equation*}\label{def-weak-sol}
	\begin{split}
	\langle u,\phi_1\rangle &+ \int_{\mb R^n}a_1(x)u|u|^{p-2}\phi_1~dx+ \langle v,\phi_2\rangle + \int_{\mb R^n}a_2(x)v|v|^{p-2}\phi_2~dx\\
	&   -\alpha \int_{\mb R^n}(|x|^{-\mu}*|u|^q)u|u|^{q-2}\phi_1 ~dx-\gamma \int_{\mb R^n}(|x|^{-\mu}*|v|^q)v|v|^{q-2}\phi_2 ~dx\\
	&  -\beta \int_{\mb R^n}(|x|^{-\mu}*|v|^q)u|u|^{q-2}\phi_1~ dx-\beta \int_{\mb R^n}(|x|^{-\mu}*|u|^q)v|v|^{q-2}\phi_2 ~dx\\
	& - \int_{\mb R^n}(f_1\phi_1 +f_2\phi_2)~dx=0,\; \forall \;(\phi_1,\phi_2) \in Y.
	\end{split}
	\end{equation*}
\end{Definition}
Thus we define the energy functional corresponding to $(P)$ as
\begin{equation*}
\begin{split}
J(u,v) &= \frac{1}{p}\|(u,v)\|^p - \frac{1}{2q}\int_{\mb R^n}\left(\alpha(|x|^{-\mu}*|u|^q)|u|^q +\beta (|x|^{-\mu}*|u|^q)|v|^q \right)dx\\
&\quad - \frac{1}{2q}\int_{\mb R^n}\left( \beta (|x|^{-\mu}*|v|^q)|u|^q+ \gamma (|x|^{-\mu}*|v|^q)|v|^q \right)dx
-\int_{\mb R^n } (f_1u+f_2v)dx\\
&= \frac{1}{p}\|(u,v)\|^p - \frac{1}{2q}\int_{\mb R^n}\left(\alpha(|x|^{-\mu}*|u|^q)|u|^q +2\beta (|x|^{-\mu}*|u|^q)|v|^q\right.\\
& \quad \left.+ \gamma (|x|^{-\mu}*|v|^q)|v|^q \right)dx -\int_{\mb R^n } (f_1u+f_2v)dx.
\end{split}
\end{equation*}
Clearly, weak solutions to $(P)$ corresponds to the critical points of $J$. To find the critical points of $J$, we constraint our functional $J$ on the Nehari manifold
\[\mc N = \{(u,v)\in Y: \; (J^\prime(u,v),(u,v)) =0 \},\]
where
\begin{equation*}
\begin{split}
(J^\prime(u,v),(u,v)) =& \|(u,v)\|^p - \int_{\mb R^n}\left(\alpha(|x|^{-\mu}*|u|^q)|u|^q +2\beta (|x|^{-\mu}*|u|^q)|v|^q\right.\\
&\left. + \gamma (|x|^{-\mu}*|v|^q)|v|^q \right)dx
-\int_{\mb R^n } (f_1u+f_2v)dx.
\end{split}
\end{equation*}
Clearly, every nontrivial weak solution to $(P)$ belongs to $\mc N$. Denote $I(u,v)= (J^\prime (u,v),(u,v))$ and subdivide the set $\mc N$ into three sets as follows:
\begin{equation*}
{\mc N^\pm = \{(u,v)\in \mc N: \; \pm(I^\prime (u,v),(u,v))> 0\}},
\end{equation*}
\begin{equation*}
\mc N^0 = \{(u,v)\in \mc N: \; (I^\prime (u,v),(u,v))=0\}.
\end{equation*}
Then $\mc N^0$ contains the element $(0,0)$ and $\mc N^+ \cup \mc N^0$ and $\mc N^- \cup \mc N^0$ are closed subsets of $Y$. For $(u,v)\in Y$, we define the fibering map $\varphi: (0,\infty) \to \mb R $ as $\varphi(t) = J(tu,tv)$. One can easily check that $(tu,tv)\in \mc N$ if and only if $\varphi^\prime(t)=0$, for $t>0$
and  $\mc N^+$, $\mc N^-$ and $\mc N^0$ can also be written as
\begin{equation*}
\mc N^\pm = \{ (tu,tv)\in \mc N:\; \varphi^{\prime\prime}(t)\gtrless 0 \},\\
\text{ and }\mc N^0 = \{ (tu,tv)\in \mc N:\; \varphi^{\prime\prime}(t)=0 \}.
\end{equation*}
We showed that $J$ becomes coercive and bounded from below on $\mc N$. By analyzing  the fiber maps $\varphi_{u,v}(t)$ we proved that if \eqref{star0} holds, then $\mc N_0 =\{(0,0)\}$ and $\mc N^-$ is a closed set. By Lagrange multiplier method, we showed that minimizers of $J$ over $\mc N^+$ and $\mc N^-$ are the weak solutions of $(P)$. So our problem reduced to minimization problem given below-
\begin{equation*}\label{min-N+}
\Upsilon^+ := \inf\limits_{(u,v) \in \mc N^+}J(u,v),\; \text{and}\;
\Upsilon^- := \inf\limits_{(u,v) \in \mc N^-}J(u,v).
\end{equation*}
Using again the map $\varphi_{u,v}$, we could show that $\Upsilon^+<0$ whereas $\Upsilon^->0$. Our next task was to consider
\[\Upsilon := \inf_{(u,v)\in \mc N}J(u,v)\]
and show that there exist a constant $C_1>0$ such that $\Upsilon \leq - \frac{(2q-p)(2qp-2q-p)}{4pq^2}C_1. $
Our next result was crucial one, which concerns another minimization problem.
\begin{Lemma}\label{inf-achvd}
	For $0\neq f_1,f_2 \in L^{\frac{p}{p-1}}(\mb R^n)$,
	\[\inf_Q \left( C_{p,q}\|(u,v)\|^{\frac{p(2q-1)}{2q-p}}- \int_{\mb R^n}(f_1u+f_2v)~dx\right):= \delta\]
	is achieved, where $ Q = \{(u,v)\in Y : L(u,v)=1\}$. Also if $f_1,f_2$ satisfies \eqref{star0}, then $\delta >0$.
\end{Lemma}
After this, using the Ekeland variational principle we proved the existence of a Palais Smale sequence for $J$ at the levels $\Upsilon$ and $\Upsilon^-$. Keeping this altogether, we could prove that $\Upsilon$ and $\Upsilon^-$ are achieved by some functions $(u_0,v_0)$ and $(u_1,v_1)$ where $(u_0,v_0)$  lies in $\mc N^+$ and forms a local minimum of $J$. The non negativity of $(u_i,v_i)$ for $i=0,1$ was showed using the modulus function $(|u_i|,|v_i|)$ and their corresponding fiber maps. Hence we conclude our main result, Theorem \ref{mainthrm}.
\subsection{Doubly nonlocal system with critical nonlinearity}
In this section, we illustrate our results concerning a system of Choquard equation with Hardy-Littlewood-Sobolev critical nonlinearity which involves the fractional Laplacian. Precisely, we consider the following problem in \cite{TS-5}
$$
(P_{\la,\delta})\left\{
\begin{array}{rlll}
(-\De)^su &= \la |u|^{q-2}u + \left(\int_{\Om}\frac{|v(y)|^{2^*_\mu}}{|x-y|^\mu}~\mathrm{d}y\right) |u|^{2^*_\mu-2}u\; \text{in}\; \Om\\
(-\De)^sv &= \delta |v|^{q-2}v + \left(\int_{\Om}\frac{|u(y)|^{2^*_\mu}}{|x-y|^\mu}~\mathrm{d}y \right)  |v|^{2^*_\mu-2}v \; \text{in}\; \Om\\
u &=v=0\; \text{in}\; \mb R^n\setminus\Omega,
\end{array} \right.$$
where $\Om$ is a smooth bounded domain in $\mb R^n$, $n >2s$, $s \in (0,1)$, $\mu \in (0,n)$, $2^*_\mu = \displaystyle\frac{2n-\mu}{n-2s}$ is the upper critical exponent in the Hardy-Littlewood-Sobolev inequality, $1<q<2$, $\la,\delta >0$ are real parameters. As we illustrated some literature on system of elliptic equation involving fractional Laplacian in the last subsection, it was an open question regarding the existence and multiplicity result for system of Choquard equation with Hardy-Littlewood-Sobolev critical nonlinearity, even in the local case $s=1$.
Using the Nehari manifold technique, we prove the following main result-
\begin{Theorem}\label{MT}
	Assume $1<q<2$ and $0<\mu<n$ then {there exists positive constants} $\Theta$ and $\Theta_0$ such that
	\begin{enumerate}
		\item if $\mu\leq 4s$ and $ 0< \la^{\frac{2}{2-q}}+ \delta^{\frac{2}{2-q}}< \Theta$, the system $(P_{\la,\delta})$ admits at least two nontrivial solutions,
		\item if $\mu> 4s$ and $ 0< \la^{\frac{2}{2-q}}+ \delta^{\frac{2}{2-q}}< \Theta_0$, the system $(P_{\la,\delta})$ admits at least two nontrivial solutions.
	\end{enumerate}
	Moreover, there exists a positive solution for $(P_{\la,\delta})$.
\end{Theorem}
Consider the product space $Y:= X_0\times X_0$ endowed with the norm $\|(u,v)\|^2:=\|u\|^2+\|v\|^2$. For notational convenience, if $u, v \in X_0$ we set
\[ B(u,v):= \int_{\Om}(|x|^{-\mu}\ast|u|^{2^*_\mu})|v|^{2^*_\mu}. \]
\begin{Definition}
	We say that $(u,v)\in Y$ is a weak solution to $(P_{\la,\delta})$ if for every $(\phi,\psi)\in Y$, it satisfies
	\begin{equation*}
	\begin{split}\label{weak-sol}
	(\langle u, \phi\rangle + \langle v,\psi\rangle) &= \int_{\Om}(\la |u|^{q-2}u\phi+\delta |v|^{q-2}v\psi)\mathrm{d}x\\
	&+ \int_{\Om}(|x|^{-\mu}\ast|v|^{2^*_\mu})|u|^{2^*_\mu-2}u\phi~\mathrm{d}x +  \int_{\Om}(|x|^{-\mu}\ast|u|^{2^*_\mu})|v|^{2^*_\mu-2}v\psi~\mathrm{d}x.
	\end{split}
	\end{equation*}
\end{Definition}
Equivalently, if we define the functional $I_{\la,\delta}:Y\to \mb R$ as
\[I_{\la,\delta}(u):= \frac{1 }{2} \|(u,v)\|^2-\frac{1}{q}\int_{\Om}(\la |u|^q+\delta |v|^q)-\frac{2}{22^*_\mu}B(u,v) \]
then the critical points of $I_{\la,\delta}$ correspond to the weak solutions of $(P_{\la,\delta})$. We set
\[\tilde S_s^H = \inf_{(u,v) \in Y\setminus \{(0,0)\}} \frac{\|(u,v)\|^2 }{ \left( \int_{\Om}(|x|^{-\mu}\ast |u|^{2^*_\mu})|v|^{2^*_\mu}~\mathrm{d}x\right)^{\frac{1}{2^*_\mu}}} = \inf_{(u,v) \in Y\setminus \{(0,0)\}} \frac{\|(u,v)\|^2}{B(u,v)^{\frac{1}{2^*_\mu}}} \]
and show that $\tilde S_s^H = 2 S_s^H$.
We define the set
\[\mc N_{\la,\delta}:= \{(u,v)\in Y\setminus \{0\}:\; (I_{\la,\delta}^\prime(u,v),(u,v))=0\}\]
and find that the functional $I_{\la,\delta}$ is coercive and  bounded below on $\mc N_{\la,\delta}$. Consider the fibering map $\varphi_{u,v}:\mb R^+ \to \mb R$ as
$\varphi_{u,v}(t)= I_{\la,\delta}(tu,tv)$
which gives another characterization of $\mc N_{\la,\delta}$ as follows
\[\mc N_{\la,\delta}=\{(tu,tv)\in Y\setminus\{(0,0)\}:\; \varphi_{u,v}^\prime (t)=0\}
\]
because $\varphi_{u,v}^\prime(t)= (I_{\la,\delta}^\prime(tu,tv),(u,v))$. Naturally, our next step is to divide $\mc N_{\la,\delta}$ into three subsets corresponding to local minima, local maxima and point of inflexion of $\varphi_{u,v}$ namely
\begin{align*}
\mc N_{\la,\delta}^\pm := \{(u,v)\in \mc N_{\la,\delta}:\; \varphi_{u,v}^{\prime\prime}(1)\gtrless 0\}\;\; \text{and}\;\;\mc N_{\la,\delta}^0 := \{(u,v)\in \mc N_{\la.\delta}:\;  \varphi_{u,v}^{\prime\prime}(1)=0\}.
\end{align*}
As earlier, the minimizers of $I_{\la,\delta}$ on $\mc N_{\la,\delta}^+$ and $\mc N_{\la,\delta}^-$ forms nontrivial weak solutions of $(P_{\la,\delta})$. Then we found a threshold on the range of $\la$ and $\delta$ so that $\mc N_{\la,\delta}$ forms a manifold. Precisely we proved-
\begin{Lemma}\label{Theta-def-lem}
	For every $(u,v)\in Y \setminus \{(0,0)\}$ and $\la,\delta $ satisfying $0<\la^{\frac{2}{2-q}}+ \delta^{\frac{2}{2-q}} < \Theta$,
	where $\Theta$ is equal to
	\begin{equation}\label{Theta-def}
	\left[ \frac{2^{2^*_\mu-1}{(C_s^n)^{\frac{22^*_\mu-q}{2-q}}}}{C(n,\mu)}\left(\frac{2-q}{22^*_\mu-q}\right) \left( \frac{22^*_\mu-2}{22^*_\mu-q}\right)^{\frac{22^*_\mu-2}{2-q}} S_s^{\frac{q(2^*_\mu-1)}{2-q}+2^*_\mu}|\Om|^{-\frac{(2^*_s-q)(22^*_\mu-2)}{2^*_s(2-q)}}\right]^{\frac{1}{2^*_\mu-1}}
	\end{equation}
	{then there exist unique} $t_1,t_2>0$ such that $t_1<t_{\text{max}}(u,v)<t_2$, $(t_1u,t_1v) \in \mc N_{\la,\delta}^+$ and $(t_2u,t_2v)\in \mc N_{\la,\delta}^-$. Moreover, $\mc N_{\la,\delta}^0= \emptyset$. As a consequence, we infer that for any $\la, \delta$ satisfying $0 < \la^{\frac{2}{2-q}} + \delta^{\frac{2}{2-q} }< \Theta$,
	\[\mc N_{\la,\delta}= \mc N_{\la,\delta}^+ \cup \mc N_{\la,\delta}^-. \]
\end{Lemma}
After this, we prove that any Palais Smale sequence $\{(u_k,v_k)\}$ for $I_{\la,\delta}$ must be bounded in $Y$ and its  weak limit  forms a weak solution of $(P_{\la,\delta})$.
We define the following
\[l_{\la,\delta}= \inf_{\mc N_{\la,\delta}} I_{\la,\delta} \; \text{and} \;l_{\la,\delta}^\pm = \inf_{\mc N_{\la,\delta}^\pm} I_{\la,\delta}. \]
We fix $0 < \la^{\frac{2}{2-q}} + \delta^{\frac{2}{2-q}}< \Theta$ and showed that $l_{\la,\delta} \leq l_{\la,\delta}^+<0$ and $\inf \{\|(u,v)\|:\; (u,v)\in \mc N_{\la,\delta}^- \}>0$.\\
To prove the existence of first solution, we first show that there exists a $(PS)_{l_{\la,\delta}}$ sequence $\{(u_k,v_k)\} \subset \mc N_{\la,\delta}$ for $I_{\la,\delta}$ using the Ekeland variational principle and then prove that $l_{\la,\delta}^+$ is achieved by some function $(u_1,v_1) \in \mc N_{\la,\delta}^+$. Moreover $u_1,v_1>0$ in $\Om$ and for each compact subset $K$ of $\Om$, there exists a $m_K>0$ such that $u_1,v_1 \geq m_K$ in $K$. Thus we obtain a positive weak solution $(u_1,v_1)$ of $(P_{\la,\delta})$.\\
On the other hand, proof of existence of second solution has been divided into two parts- $\mu\leq 4s$ and $\mu>4s$. In the case $\mu \leq 4s$, using the estimates in Proposition \ref{estimates1}, we could reach the first critical level as follows-
\[ \sup_{t \geq 0} I_{\la,\delta}((u_1,v_1)+ t(w_0,z_0) )< c_1:=I_{\la,\delta}(u_1,v_1)+ \frac{n-\mu+2s}{2n-\mu} \left(\frac{{C_s^n} \tilde S_s^H}{2} \right)^{\frac{2n-\mu}{n-\mu+2s}} \]
for some non negative $(w_0,z_0) \in Y\setminus \{(0,0)\}$. This implied $l_{\la,\delta}^- < c_1$. Whereas to show the same thing in the case $\mu>4s$, we had to take another constant $\Theta_0 \leq \Theta$ and the same estimates as in Proposition \ref{estimates1}. Consequently, we prove that there exists a $(u_2,v_2) \in \mc N_{\la,\delta}^-$ such that $l_{\la,\delta}^-$ is achieved, hence gave us the second solution. From this, we concluded the proof of Theorem \ref{MT}.

\section{Some open questions}
Here we state some open problems in this direction.
\begin{enumerate}
	\item $H^1$ versus $C^1$ local minimizers and global multiplicity result: Consider  energy functional defined on $H^1_0(\Om)$ given by $\Phi(u)=\frac{\|u\|^2}{2}- \la \int_\Om F(x,u)$ where $F$ is the primitive of $f$. When $|f(u)|\leq C(1+|u|^{p})$ for $p \in
	[1,2^*]$, Brezis and Nirenberg in \cite{brenir} showed that a local minimum of $\Phi$ in $C^1(\Om)$-topology is also a local minimum in the  $H^1_0(\Om)$-topology. Such
	property of the functional corresponding to Choquard type nonlinearity and singular terms is still not addressed.
	\item Variable exponent problems:  As pointed out in section \ref{sec-1.2}, existence of a solution for problem \eqref{BS-12} has been studied in \cite{p(x)-choq}
	but the question of multiplicity of solutions for variable exponent Choquard equations is still open.
	\item $p$-Laplacian critical problems: The Critical exponent problem involving the $p$-Laplcian and Choquard terms is an important question. This requires the study of minimizers of $S_{H,L}$. Also the regularity of solutions and the global multiplicity results for convex-concave nonlinearities is worth exploring.
	\item Hardy-Sobolev operators and nonlocal problems: The doubly critical problems arise due to the presence of two noncompact terms. Hardy Sobolev operator is defined as $-\De_p u -\frac{\mu |u|^{p-2}u}{|x|^2}$. Here the critical growth Choquard terms in the {equations} require the minimizers and {asymptotic} estimates to study the compactness of minimizing sequences. The existence and multiplicity results are good questions to explore in this case.
\end{enumerate}

\end{document}